\documentclass{birkjour}

\usepackage{amsfonts,amssymb,amscd, amsmath}
\usepackage{epsfig,psfrag}  
\usepackage{graphicx}

\newtheorem{teo}{Theorem}[section]
\newtheorem{prop}[teo]{Proposition}
\newtheorem{lema}[teo]{Lemma}

\newtheorem{coro}[teo]{Corollary}
\newcommand{\C}{{\mathbb C}}
\newcommand{\R}{{\mathbb R}}

\newcommand{\Z}{{\mathbb Z}}
\newcommand{\N}{{\mathbb N}}

\newcommand{\fol}{{\mathcal F}}
\newcommand{\tilf}{{\widetilde{\mathcal F}}}

\newcommand{\mathcalp}{{\mathcal P}}

\newcommand{\mathcals}{{\mathcal S}}

\newcommand{\palix}{{\partial / \partial x}}
\newcommand{\paliy}{{\partial /\partial y}}
\newcommand{\mcd}{\mathcal{D}}

\begin{document}

\received{January 11, 2011}
\revised{May 25, 2011}
\accepted{June 01, 2011}

%
%
%
%
%
%
%
%
%

\title[Fibrations of genus two on complex surfaces]
 {Fibrations of genus two on complex surfaces}

\author[Julio Rebelo]{Julio Rebelo}

\address{%
1- Universit\'e de Toulouse, UPS, INSA, UT1, UT2 \\
Institut de Math\'ematiques de Toulouse\\
F-31062 Toulouse, FRANCE.\\
2- CNRS, Institut de Math\'ematiques de Toulouse UMR 5219\\
F-31062 Toulouse, FRANCE.\\}

\email{rebelo@math.univ-toulouse.fr}

\thanks{The second author was partially supported by the NSF grant DMS 1007155.}
\author{Bianca Santoro}
\address{Mathematics Department \\
The City College of New York \\
138th St. and Convent Avenue, NAC 6/203C\\
New York, NY 10032\\
USA\\
}
\email{bsantoro@ccny.cuny.edu}


\date{May 27th, 2011.}

\dedicatory{Dedicated to Keti Tenenblat on the occasion of her $65^{\rm th}$ birthday.}

\begin{abstract}
This work is devoted to the study of fibrations of genus~$2$ by using as its main tool the
theory of singular holomorphic foliations. In particular we obtain a sharp differentiable version of Matsumoto-Montesinos
theory. In the case of isotrivial fibrations, these methods are powerful enough to provide a detailed global
picture of the both the ambient surface and of the structure of the fibrations itself.
\end{abstract}


\maketitle
\section{Introduction}

A complex surface $M$ is called ruled (resp. elliptic) if it
carries a singular fibration $M \rightarrow S$ of genus~$0$ (resp.~$1$, cf. Section~2). The
structure of ruled surfaces was known to the classical italian
geometers whereas Kodaira has provided a similar picture for elliptic surfaces. After Kodaira, several authors have considered the
case of higher genus fibrations. Among these works, we can distinguish between papers concerned with the global geometry of the ambient
surface $M$ carrying a particular fibration (for example the existence of special relations among its invariants) and papers that focus on
the structure of the pencil itself. The present work certainly belongs to the latter category as it will become clear below.
In fact a standard point of view consists of splitting the discussion in two parts as follows:
\begin{enumerate}
\item (The local on the target point of view) Description of the fibration on a ``tubular neighborhood of a given fiber''.

\item (Globalization) The description of how the ``models'' mentioned above can be patched together in a global fibration.
\end{enumerate}
Item~1 above amounts to understanding ``fiber germs'' or ``germs of fibrations'' in the terminology of \cite{reid}. Concerning fibrations
of genus~$2$, Ogg \cite{ogg}, and independently Iitaka, has provided a list of all possible singular fibers. Later 
Namikawa and Ueno announced in \cite{Ueno1} the classification of the ``structure'' of these fibers. They define three
``discrete'' (explicitly computed) invariants associated to a local fibration and allowing to recognize a singular fiber unambiguously
though the ``same'' fiber may appear in families with different invariants. In a sense these results can be regarded as the ``infinitesimal''
version of problem~1 above: whereas they describe the structure of a singular fiber, they fall short of characterizing the fibration on a
tubular neighborhood of it. For readers more familiar with algebraic methods, the distinction is somehow similar to the difference between ``stalk'' and
``fiber'' of a sheaf, as pointed out in \cite{reid}.

The topological version of problem~1 was solved by Matsumoto and Montesinos.
They showed that the topological type of a neighborhood of a singular fiber in a fibration of genus~$\textsf{g}$
($\textsf{g} \geq 2$) is determined by the action of the monodromy in the mapping class group of the surface $S_{\textsf{g}}$ of genus~$\textsf{g}$. In other words,
the topological type is determined by a sort of ``non-commutative monodromy'' where $H_1 (S_{\textsf{g}} ; \Z)$ is substituted by the
fundamental group $\pi_1 (S_{\textsf{g}})$ of $S_{\textsf{g}}$
(recall also that the mapping class group of $S_{\textsf{g}}$ may be identified to
${\rm Aut}\, (\pi_1 (S_{\textsf{g}}))/ {\rm Inn}\, (\pi_1 (S_{\textsf{g}}))$). The context here is as follows: consider two ``local'' fibrations of genus~$2$ (or higher)
over the disk $D \subset \C$, namely
$\mathcalp_1 : M_1 \rightarrow D$ and $\mathcalp_2 : M_2 \rightarrow D$,  whose unique singular fibers are those
sitting over $0 \in D$. A homeomorphism $H$ (resp. $C^{\infty}$-diffeomorphism) is said to {\it conjugate}\, the fibrations $\mathcalp_1, \, \mathcalp_2$ if it
preserves the structure of the fibration, i.e. if it sends fibers of $\mathcalp_1$ onto fibers of $\mathcalp_2$ (and $\mathcalp_1^{-1} (0)$ onto
$\mathcalp_2^{-1} (0)$). The two fibrations $\mathcalp_1, \, \mathcalp_2$ are then called topologically conjugate (resp. $C^{\infty}$-conjugate). The theory
of Matsumoto-Montesinos can be summarized by saying that $\mathcalp_1, \, \mathcalp_2$ are topologically conjugate if and only if their ``non-comutative
monodromy maps'' coincide. Their method also allows one to to recover the classification of \cite{Ueno1} and, in particular,
the fact that all the corresponding models actually exist, a statement originally due to Winters. In \cite{konno} the reader
will find a particularly nice account of Matsumoto-Montesinos theory along with a good deal of information and questions about fibrations
(from the point of view that focus more on the structure of the ambient manifold $M$ than on the structure of the fibration itself).
Roughly speaking, the purpose of this article may be described as an extension of Matsumoto-Montesinos results in the sense that
we shall be concerned with the regularity properties
of the conjugating homeomorphism. In fact their techniques yield little information on how regular the conjugating homeomorphism can be made.
For example it is unclear if/when the homeomorphism is actually a $C^{\infty}$-diffeomorphism (or even a holomorphic diffeomorphism).
This question is clearly important for the computation of invariants sensitive to the differentiable structure of $M$. In these
directions our main result is Theorem~A below which provides a sharp statement concerning this type of question.

To state our main result, the {\it combinatorial data of a singular fiber}\, will consist of
the structure of its irreducible components, along with their
intersections and self-intersections numbers. This data might also be called the ``decorated dual graph'' of the
fiber in question. All possible combinatorial data of fibers of a genus~$2$ fibrations were first classified in \cite{ogg}
and we shall return to this point later. Next let
$\mathcalp_1 : M_1 \rightarrow D$ and $\mathcalp_2 : M_2 \rightarrow D$ be two ``local'' fibrations of genus~$2$ as above.
The first natural question is to decide whether or not the combinatorial data determines the fibered neighborhood, i.e. whether or not the coincidence of the
combinatorial data of $\mathcalp_1^{-1} (0), \, \mathcalp_2^{-1} (0)$ guarantees the existence of a conjugating homeomorphism as above. Similarly there is the
the question of deciding whether or not a {\it same}\, fiber my appear in different degenerating families. In the case of genus~$1$, and if multiple fibers
are excluded, then the analogous question has an affirmative answer. For genus~$2$ the answer is however known to be negative as pointed out
in \cite{Ueno1}. For example in the list given in \cite{Ueno1} the model named $[2\!-\!2I_0\!-\!m]$ (page 159) and the model $[3\!-\! II_{n-0} \; \ast]$ (page 172) have identical
singular fibers whereas there is no homeomorphism conjugating their corresponding fibrations.
The nature of the latter phenomenon will be understood in the course of our discussion. However to do this, we shall be led to introduce
some invariants related to the {\it holonomy of cycles}. Let us first divide the singular fibers
in three categories as follows. Recall that fibers my have one or two irreducible components that are elliptic curves. These elliptic curves may also appear
in a singular fashion, namely as a ``collar of rational curves'' or as ``pinched torus'' these forms will be called {\it degenerate elliptic components},
cf. also \cite{ogg} or Section~2.2. Next note that a regular
elliptic component possesses two independent generators in its fundamental group. Each of these generators will yield
a positive integer number corresponding to the order of the holonomy map it gives rise. In the case of a degenerate elliptic component, it will be seen that
one of these holonomy maps is necessarily trivial so that we have only one invariant (rather than two) to keep track.
The set formed by these integers is therefore constituted by (at most) four positive integers and it will be referred to as the {\it elliptic holonomy
of the singular fiber}\, (or simply as elliptic holonomy when no misunderstood is possible). Alternatively, we can think of the singular fiber as a
``reducible singular curve'' (ie. an arrangement of possibly singular curves without taking into account any multiplicity or local index). The fundamental
group of this ``reducible singular curve'' may have up to four generators and these account for the elliptic holonomy invariant mentioned above. For example a
collar of rational curves has only one homotopically non-trivial loop, hence only one elliptic holonomy invariant will emerge (as opposed to the two generators
that occur at a regular elliptic component). The same phenomenon occurs for a pinched torus.
Our first result will be the proof that the elliptic holonomy constitutes the complementary information
needed to explain the appearance of a ``same singular fiber'' in different fibrations mentioned above. In particular we shall prove the following:

\vspace{0.2cm}

\noindent {\bf Theorem A}. {\sl Let ${\mathcal P}_1 : M_1 \rightarrow D$ and ${\mathcal P}_2 : M_2 \rightarrow D$
be two fibrations of genus $2$ over the disk $D \subset \C$ whose unique singular fibers are those
sitting over $0 \in D$. Suppose that
the combinatorial data of the singular fiber ${\mathcal P}_1^{-1}(0)$ is the same as 
${\mathcal P}_2^{-1}(0)$ and that all their irreducible components are rational curves. Suppose also that the corresponding Dynkin diagram contains no loop.
Then there exists a $C^{\infty}$-diffeomorphism that conjugates the 
two fibrations. Furthermore, this diffeomorphism is transversely holomorphic.

More generally if $\mathcalp_1, \, \mathcalp_2$ are general singular fibers with the same combinatorial data, then the above statement still holds
provided that their elliptic holonomy invariants coincide.}

\vspace{0.1cm}

\noindent {\bf Remark}. {\rm The combinatorial data cannot holomorphically determine the structure of the fibration since the complex
structure of the regular fibers may vary. Yet Proposition~\ref{localmodels} shows that
it holomorphically determines the structure of the singularities of the singular fiber whether they are thought of as
singularities of a ``reducible analytic curve'' or as singularities of a local foliation (given by the restriction of the corresponding fibration,
cf. Section~2.2 for further details). This is a typical statement that concerns the structure of the fibration itself more than the global nature
of the ambient manifold $M$. The contents of Proposition~\ref{localmodels} however become false for fibrations of sufficiently high genus.}

\vspace{0.1cm}

With respect to the statement of Theorem~A, it should be noted that the actual values that can be attained by the elliptic holonomy will be made explicit for each
fixed combinatorial data for the singular fiber, as shown at the end of Section~3. Thus the assumption made in Theorem~A is necessary and does not reduce
its applicability. A corollary of this discussion is as follows:

\vspace{0.2cm}

\noindent {\bf Complement to Theorem A}. {\sl For a fixed combinatorial data of a singular fiber, its elliptic holonomy invariant can take at most on
two values. In particular a ``same'' singular fiber can belong at most to two different families of degenerating genus~$2$ curves. Besides if the elliptic
holonomy invariant takes on more than on single value, then the singular fiber in question must contain a (possibly singular) elliptic curve as an irreducible
component or a loop of rational curves.}

\vspace{0.1cm}

Naturally the preceding results fit into the pattern of approaching fibrations by first understanding ``tubular neighborhoods of fibers''. Indeed we are
interested in describing the local structure of the fibration as accurately as possible. It is our hope that the information collected above will find further
applications in the study of complex surfaces carrying a genus~$2$ fibration. In this paper, however, the only global case that will be discussed
correspond to the special case of {\it isotrivial fibrations} i.e. fibrations whose fibers are pairwise isomorphic as Riemann surfaces. In fact,
an obvious consequence of Theorem~A is that the obstruction for the existence of a topological conjugacy between $\mathcal{P}_1, \, \mathcal{P}_2$
coincides with the obstruction for the existence of a $C^{\infty}$-conjugacy (being, in addition, transversely holomorphic). It is then natural to wonder under what
conditions the existence of a holomorphic conjugacy can be ensured. For this let us assume that both $\mathcal{P}_1, \, \mathcal{P}_2$ are isotrivial fibrations.
Clearly the complex structure of $\mathcal{P}_1^{-1} (z)$ (resp. $\mathcal{P}_1^{-1} (z)$) is an invariant in this case. Next we have:

\vspace{0.2cm}

\noindent {\bf Theorem B}. {\sl Let ${\mathcal P}_1 : M_1 \rightarrow D$ and ${\mathcal P}_2 : M_2 \rightarrow D$
be conjugate by a diffeomorphism $H$ as in Theorem~A. Suppose that $\mathcal{P}_1$ (resp. $\mathcal{P}_2$)
is isotrivial and that the corresponding generic fibers $\mathcal{P}_1^{-1} (z)$, $\mathcal{P}_2^{-1} (z)$ are isomorphic Riemann surfaces.
Then $\mathcal{P}_1, \, \mathcal{P}_2$ are conjugate by a holomorphic diffeomorphism.}

\vspace{0.1cm}

After having proved the above theorems, we shall close the paper by providing a general description of isotrivial fibrations of genus~$2$. This
description will take up most of the contents of Section~4 but, for the time being, we shall content ourselves of stating the Theorem~C below
(which is certainly not sharp). Note that the assumption that the singular fibers have only normal crossings do not affect the generality of the statement,
cf. Section~2 or Section~3.

\vspace{0.2cm}

\noindent {\bf Theorem C}. {\sl Let $\mathcalp : M \rightarrow S$ be an isotrivial fibration of genus~$2$ all of whose singular fibers have only
normal crossing singularities. Then there exists
a cyclic covering $\widetilde{S} \rightarrow S$, ramified only over the critical values of $\mathcalp$ in $S$ and over a single additional point in $S$ such that
the fiber product $\widetilde{S} \times_S M$ is bimeromorphically equivalent to a regular fibration of genus~$2$ (necessarily isotrivial).

Moreover the only prime factors that can divide the multiplicities of the singular fibers of $\mathcalp$ are~$2, \, 3$ and~$5$. The automorphism group of the
typical fiber is either cyclic or it belongs to the following list: the Alternating group $\mathcal{A}_4$, the Symmetric group $\mathcal{S}_4$ or the Dihedral
group $D_n$ for $n \in \{ 6, 8, 9, 10, 12, 15, 16, 18, 20\}$. Finally if one among the irreducible components of the singular fibers of $\mathcalp$ there is one whose
multiplicity has~$5$ as a prime factor, then the mentioned automorphism group is either cyclic (with an order that is a multiple of~$5$) or it is one of the
Dihedral groups $D_{10}, \, D_{15}, \, D_{20}$.}

\vspace{0.1cm}

From a more technimathcal Point of view, the upshot of this paper is definitely the systematic use of the standard theory of (singular) holomorphic foliations. This is a geometric
theory that also has a quantitative analytic character. In fact, in some sense, the topological arguments considered by Matsumoto-Montesinos
\cite{MMontesinos} are ``purely qualitative''. This means that their theory, while extremely effective for the description of topological models and for
monodromy computations, is not well-adapted to answer ``more quantitative'' questions involving
for example the regularity of a conjugating
homeomorphism. The advantage of our approach lies exactly in the ``analytic'' (or differentiable) character intrinsically attached to the theory of foliations
which becomes a suitable tool for conducting a finer study of the mentioned topological models, leading in particular to Theorems~A and~B.

Alerted by the referee of this paper, the authors have recently become aware of a number of papers studying fibrations of small genus, cf.
\cite{konno}, \cite{catanese}, \cite{horikawa}, \cite{reid}, \cite{xiao} where the reader will find results concerned with
the ``geography'' of manifolds carrying these fibrations as well as results and methods that are aimed at understanding the structure
of fibrations itself. Our previous experience seems to indicate that this literature is not widely known to colleagues working in (differential)
geometry and topology (and even on certain complex dynamical systems). In fact, for this type of public (including the authors of the present article) the algebraic
techniques on which most of the mentioned papers are based are harder to grasp with than the methods inspired in foliation theory that
are employed here. Hopefully this lack of communication will soon be bridged and we also hope that our techniques will bring a
contribution to some questions raised in former works. For example, it is reasonable to expect that they will have a saying concerning the ``Morsification''
of singular fibers and related problems about the (local) relatively canonical algebra according to the point of view of M. Reid and G. Xiao.

Let us close this Introduction with an outline of the contents of this paper. The general strategy for Theorem~A consists of three main steps, the first two being
the classification of the singular points of $\mathcalp$ in the context of ``local singular foliations'' and the description of the Dynkin diagram of the corresponding
singular fiber. Finally the third step consists of working out the ways in which the ``local models for $\mathcalp$'' at its singular points can be glued together
following the pattern given by a fixed Dynkin diagram.

The characterization of all Dynkin diagrams associated to the singular fibers is given in \cite{ogg} and its reviewed in Section~2.1. This review is useful not
only to make this paper more self-contained but also because the corresponding algorithms play a role in Section~3, specially in the proof of Proposition~\ref{localmodels}.
Section~2.2 begins with some well-known results related to indices of singularities of holomorphic foliations, a general reference where the reader will find
foundational material about singular foliations as well as all the results used in this paper is \cite{reference?}. The discussion then continues with the
classification of all singularities of $\mathcalp$ (viewed as a local singular foliation). Thus the first two steps mentioned above are carried out in Section~2.
Some of the material presented here is rather elementary and/or might be skipped
modulo referring the reader to previous work. Alternatively the discussion might be made slightly shorter by
building on results more elaborated than the original classification due to Ogg. We decided not to follow any of these possibilities in order to keep
the section as elementary and self-contained as possible. Hopefully this will help to make the subject of low genus fibrations more accessible
to a public with different background.

Section~3 contains the third step leading to Theorem~A. First we prove Proposition~\ref{localmodels} that roughly speaking asserts that ``corresponding
singularities of genus~$2$ fibrations lying in isomorphic Dynkin diagrams are themselves isomorphic as singularities of foliations''. To go from this statement
to Theorem~A it will require us to discuss how these singularities can be glued together by means of ``pieces of regular foliations containing the regular part
of the singular fiber''. The key notions allowing us to work out the ways in which this gluing is made are those of ``holonomy representation 
generated at a given component of the singular fiber'' and the ``total holonomy representation of the singular fiber''. All this is detailed in the course of
Section~3 that ends with the proof of Theorem~A and of its complement.

Section~4 is devoted to the proofs of Theorems~B and~C. This section begins with the proof of Theorem~B (relying on Theorem~A). The proof
of Theorem~C is given in the course of the subsequent global description of isotrivial fibrations.
The starting point of our analysis is Proposition~\ref{Lastsection1.1}, valid for
general isotrivial fibrations, stating the existence of another fibration that is ``transverse'' to $\mathcalp$ in an appropriate sense. This ``transverse''
fibration allows us, in particular, to represent standard monodromy by a holomorphic diffeomorphism. Furthermore it also allows us to bring
all these representative of monodromy maps together into a single group of automorphism of a generic fiber. The resulting structure turns out to
be very rigid and it yields a large amount of information on the fibration in question.

A last comment is that we did not include a list of all possible
types of singular fibers and discrete invariants as in Namikawa et Ueno \cite{Ueno1} since this would lead to a pointless duplication
of published work. However for the reading of this paper, and specially of Section~3, is highly recommend to have their list (or Ogg's list)
in hand.

\section{Some background material}

A singular fibration on a compact complex surface $M$ is a
non-constant holomorphic map $\mathcalp$ from $M$ to a compact
Riemann surface $S$. Given one such map $\mathcalp$, the finite
set $\{ p_1 ,\ldots, p_s \} \subset S$ consisting of the
critical values of $\mathcalp$ is such that $\mathcalp$ defines a regular
fibration of  $M \setminus \bigcup_{i=1}^s \mathcalp^{-1} (p_i)$
over $S \setminus \{ p_1 ,\ldots, p_s \}$. The {\it genus}\, of
$\mathcalp$ is simply the genus of the fiber of the mentioned (regular)
fibration or, equivalently, the genus of a ``generic'' fiber of $\mathcalp$.
The preimages $\mathcalp^{-1} (p_i)$ of the critical values of $\mathcalp$
are called {\it singular fibers}. Let us first briefly review the work of Ogg \cite{ogg} concerning the
Dynkin diagrams associated to singular fibers in genus~$2$ fibrations.

\subsection{Dynkin Diagrams - The structure of the singular fiber}

Recall that Kodaira \cite{kod} has classified all possible singular fibers for an elliptic fibration. Besides an elliptic ``multiple'' fiber, these are either a rational
curve with a node or a cusp (the so-called pinched torus), or a special arrangement of rational curves with
self-intersection $-2$.
Kodaira's argument is elementary in nature. It was extended by Ogg  \cite{ogg} and, independently, by Iitaka, to fibrations of genus~$2$.
For the convenience of the reader, and also because it fits the purpose of this paper of presenting a systematic
treatment of questions relative to singular fibrations, we are going to summarize Ogg's classification here.
To begin with, let $\mathcalp: M \rightarrow S$ denote again a fibration consisting of curves
of genus $\textsf{g}$, with connected fibers. Throughout this section, the fibers of $\mathcalp$ are supposed to be {\it minimal}\,
in the sense that they do not contain a irreducible component consisting of a rational curve with self-intersection~$-1$.
Note that this assumption does not force $M$ to be minimal as well, in other words $M$ may contain
$-1$ rational curves. Yet, this suffices to allows us to treat fibrations and pencils on the same footing.

Note that every two fibers are equivalent as divisors, i.e. they 
belong to the same homology class whether or not they are singular. Hence, the self-intersection of any
fiber $L$ is $L\cdot L=L^2 = 0$. Moreover, if $D$ belongs to the group generated by the irreducible
components of a fiber $L$, we still have $L \cdot D = 0$. Next, we recall that the
{\it arithmetic genus} of a fiber $L$ is defined by 
$$
\textsf{g} \, (L) = 1 + \frac{1}{2}(L^2 + L \cdot K),
$$
where $K$ stands for the  canonical divisor of $M$. In particular, if {\bf $L = n D $}, then the condition of
$\textsf{g}$ being an integer forces $n = 1$ in the case of $\textsf{g} = 2$.
So, if a fiber has  components of a single type, it  can only contain one single component. 
On the other hand, if a fiber $L$ is of the form {\bf $L = n D + \Gamma$}, 
$D$ and $\Gamma$ distinct, then 
$0 = L \cdot \Gamma  = n D \cdot \Gamma + \Gamma^2$. In turn, this equation implies
that $\Gamma^2 < 0$.
Finally, if {\bf $L = \sum_i n_i\Gamma_i$}, where all $\Gamma_i$'s are
distinct components, we have
$\Gamma_i \cdot K \geq 0$, and $2 \textsf{g} - 2  = \sum_i n_i \Gamma_i \cdot K$.

Now, by exploiting the fact that the arithmetic genus $\textsf{g}\, (\Gamma_i)$ of every irreducible component $\Gamma_i$
must be a non-negative integer, we conclude that $\Gamma_i^2$ is odd exactly when $\Gamma_i$ is odd.
Also, $\Gamma_i^2 \geq -2 - \Gamma_i \cdot K$.

Summarizing, setting the genus $\textsf{g}$ of the whole fiber $L$ equal to $2$, the above
relations yield the following five possibilities for the structure of a singular fiber:

\begin{description}
\item{Type A:}\hspace{1cm}
$ \Gamma_i \cdot K = 1$ \hspace{1cm} $\Gamma_i^2 = -1$  \hspace{1cm} $\textsf{g} = 1$.
\item{Type B:}\hspace{1cm}
$ \Gamma_i \cdot K = 1$ \hspace{1cm} $\Gamma_i^2 = -3$  \hspace{1cm} $\textsf{g} = 0$.
\item{Type C:}\hspace{1cm}
$ \Gamma_i \cdot K = 2$ \hspace{1cm} $\Gamma_i^2 = -2$  \hspace{1cm} $\textsf{g} = 0$.
\item{Type D:}\hspace{1cm}
$ \Gamma_i \cdot K = 2$ \hspace{1cm} $\Gamma_i^2 = -4$  \hspace{1cm} $\textsf{g} = 1$.
\item{Type E:}\hspace{1cm}
$ \Gamma_i \cdot K = 0$ \hspace{1cm} $\Gamma_i^2 = -2$  \hspace{1cm} $\textsf{g} = 0$.
\end{description}

From the relation $\sum_i n_i (\Gamma_i \cdot K) = 2\textsf{g} - 2 = 2$, we see that a 
reducible fiber of a genus $2$ fibration has only one out of five choices:

\begin{description}
\item{i)}
It has a component of type $C$, and all other fibers are of type $E$.
\item{ii)}
It has a  componentof type $D$, and all other fibers are of type $E$.
\item{iii)}
It has a  component of type $A$ with multiplicity $2$, and all other fibers are of type $E$.
\item{iv)}
It has a  component of type $A + B$, and all other fibers are of type $E$.
\item{v)}
It has a  component of type $B$ with multiplicity $2$, and all other fibers are of type $E$.
\end{description}

The strategy from now on is to study separately each of the five cases above.
The possible intersection numbers of different components will give rise to  combinatorial relations that will determine the shape of the
corresponding singular fibers.

Let us start with the case i) which is the simplest one. In other words, $L$ contains a component of type $C$ with 
multiplicity $1$, and all its other components are of type $E$. Here we first note that
$\Gamma \cdot (L - \Gamma) = 2$, so that $\Gamma$ intersects the
rest of the fiber twice. Hence, by the argument of Kodaira,
it can only be one of the seven types of cycles (of self-intersection $-2$)
described in \cite{kod}, where one of the components of multiplicity $1$ 
is replaced by $\Gamma$. This concludes the classification of
case i).

The reader should note that, in what follows, we shall keep the numbering used by Ogg in \cite{ogg} when
referring to singular fibers. In particular, the example above is  said to be of ``Type 1'' in \cite{ogg}.

\begin{figure}[hbtp]
\centering

\psfrag{G}{$\Gamma$}
\psfrag{1}{$1$}

\includegraphics[scale=0.4]{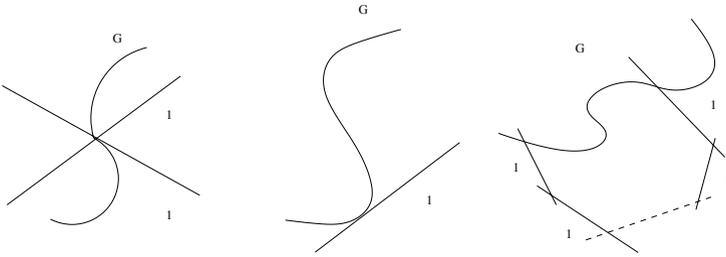}
\caption{Type $1$ singularities}
\label{Kodaira2}
\end{figure}

Let us continue by analyzing the case ii). Therefore $L$ is supposed to contain a component $\Gamma$ of type $D$.
Then, $\Gamma \cdot (L - \Gamma) = 4$, i.e. 
$\Gamma$ intersects the rest of the fiber four times more.
Certainly, one possibility for $L$ is to have $\Gamma$ 
joining two ``Kodaira components'' (Type 2), which  will be excluded henceforth.

Assume that $\Gamma_1$ (of type E) intersects $\Gamma$. We claim that
$\Gamma_1 \cdot \Gamma \leq 2$, for $0 = L \cdot \Gamma_1 = \Gamma_1^2 + \Gamma_1 \cdot \Gamma + \Gamma_1 \cdot \Gamma_2$,
and $L$ is connected (so $\Gamma_1 \cdot \Gamma_2 \geq 0$). However,
if {\bf $\Gamma_1 \cdot \Gamma = 2$}, then either 
$\Gamma_1$ intersects $\Gamma$ in two points, or the must have a double contact.
If $\Gamma_1$ appears with multiplicity $1$, then the fiber ends at $\Gamma_1$, 
and we must be in the cases already counted as Types 1 and 2.
If $\Gamma_1$ has multiplicity $2$, then $\Gamma_1$ must meet $2 \Gamma_2$, and
so on, giving rise to Type 3.
Therefore we may assume that $\Gamma$ (of Type D) intersects all the other components in at most $1$ point.

Consider now the case $\Gamma_1 \cdot \Gamma = 1$. Let $m$ be the multiplicity of $\Gamma_1$ in $L$,
$m \in \{ 1,2,3,4\}$. Then, $\Gamma_1$ meets the rest of the fiber in 
$\Gamma_i \cdot(L - m \Gamma_i) = 2m -1$ more times. Under these conditions, we have the following
algorithm to determine the possible models, i.e. the possible multiplicities
of the irreducible components:

\noindent
{\bf Step 0:} Set $\ell^0 = m$, $\ell^{-1} = 1 =$ (number of times that 
$\Gamma_1$ meets $\Gamma$).

\noindent
{\bf Step 1:} Write $2\ell^0 -1$ (the number of times $\Gamma_1$ intersects the 
rest of the fiber) as a sum of positive integers
$$
2\ell^1_1 + \dots + \ell^1_{k_1},
$$
such that $2\ell^1_j \geq \ell^0$ for all $j \in \{ 1, \cdots, k_1\}$.

If for all $j \in \{  1, \dots, k_1\}$ we have that
$2\ell^1_j = \ell^0$, then the fiber is complete, and we can stop here.
Otherwise, there exists at least one $\ell^1_{k_1, i_1}$ such that 
$2\ell^1_{k_1, i_1} - \ell^0 >0$.

\noindent
{\bf Step 2:}
For all $\ell^1_{k_1, i_j}$ from Step 2, we will write
$$
2 \ell^1_{k_1, i_j} - \ell^0 = \ell^2_{1,i_j} + \dots + \ell^2_{k_2,i_j},
$$
such that $2 \ell^2_{p,i_j}\geq \ell^1_{k_1, i_j}$, for all $p \in \{ 1, \cdots, k_2\}$.

We stop if equality holds for all $p$. Otherwise, we proceed inductively to

\noindent
{\bf Step s:}
We write, for  all $\ell^{s-1}_{k_{s-1}, i_j}$ from Step s - 1,
$$
2 \ell^{s-1}_{k_{s-1}, i_j} - \ell^{s-2}_{k_{s-2}, i_j} = \ell^s_{1,i_j} + \dots + \ell^s_{k_s,i_j},
$$
such that $2 \ell^s_{p,i_j}\geq \ell^{s-1}_{k_{s-1}, i_j}$, for all $p \in \{ 1, \cdots, k_{s}\}$.

There are three possibilities for this algorithm. Either it ends in finite steps, generating 
explicitly the Types 3 until 11 in \cite{ogg} or it generates an obvious impossibility at finite 
time, case that should be discarded. Finally a last possibility, for $s$ large enough, 
is to have the trivial decomposition of integers. This would represent 
the fiber with an infinite number
of components of type $E$, with multiplicities $m, 2m-1, 2(2m-1) - m, \dots$, what is
also impossible.

To illustrate the method above, let us study, for example, the case when 
$m = 4$.
We can write $8 - 1 = 7$ as $7$, $3 + 4$,  $2 + 5$, or $2 + 2 + 3$.
The latter three, if we keep running the algorithm, will generate the following three
models:
\begin{figure}[hbtp]
\centering
\psfrag{1}{$1$}
\psfrag{2}{$2$}
\psfrag{3}{$3$}
\psfrag{D}{$\Delta$}
\psfrag{4}{$4$}
\psfrag{5}{$5$}
\psfrag{6}{$6$}
\includegraphics[scale=0.4]{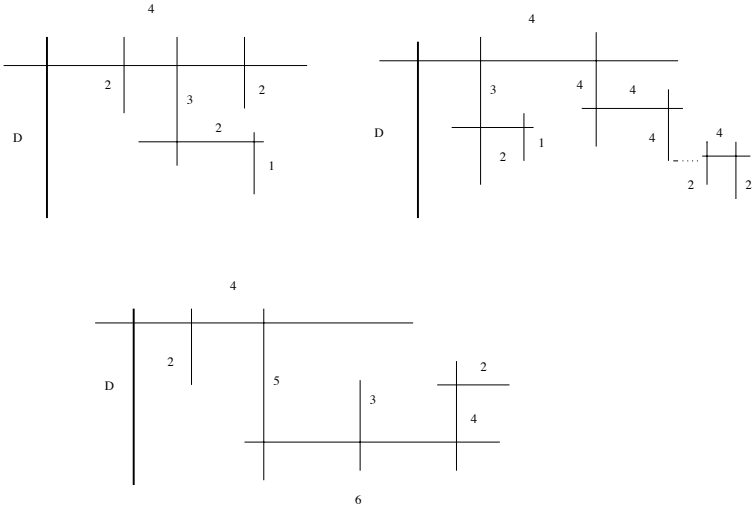}
\caption{If a fiber contains a component of type $D$ meeting a component of multiplicity $4$}
\label{3ModelsOgg}
\end{figure}

For the trivial decomposition $7 = 7$, $\Gamma_1$ will meet the fiber 
$14 - 4 = 10$ more times.
We can write $10$ as $5 + 5$, $4 + 6$, or $10$.
The former ones lead to impossibilities. Hence, 
$\Gamma_1 $ meets $10 \Gamma_2$, $\Gamma_2$ also of
type $E$.
Finally writing $20 - 7 = 13$ as $5+ 8$, $6+7$, or $13$ will clearly lead to
Type 7, an impossibility, or an infinite chain.

The cases for the other values of $m$ being completely analogous, the discussion of singular fibers
of type ii) is over.

The fibers that only contains components of Types A, B and E can be
described algorithms similar to the one presented above. It is only needed
to make further considerations with respect to the number of different components
that a fixed one meets in the singular fiber. 
Explicitly describing those algorithms, though, would make the notation unnecessarily 
heavy, and we prefer to omit it. For the complete list of possible examples, see \cite{ogg}.

\subsection{Singularities of foliations and related indices}
\label{ss.SingsandIndices}

Consider again a fibration $\mathcalp : M \rightarrow S$ with genus~$2$. The critical values of $\mathcalp$
form a finite set $\{ p_1 ,\ldots, p_s \}$.
Fixed a critical value $p_i \in S$ and a neighborhood $V_i
\subseteq S$ of $p_i$, the set $\mathcalp^{-1} (V_i)$ is going to be
called a {\it tubular neighborhood}\, of $\mathcalp^{-1} (p_i)$. To
a large extent, describing the structure of a singular fibration
$\mathcalp: M \rightarrow S$ is tantamount to providing models for
``small'' tubular neighborhoods of its singular fibers. In fact,
one passes from the set $\bigcup_{i=1}^s \mathcalp^{-1} (V_i)$ to
the whole surface $M$ by filling the complement
$M \setminus \bigcup_{i=1}^s \mathcalp^{-1} (p_i)$ in with a regular
fibration. Naturally, by ``models'' we mean ``normal forms'' for
the restriction of the fibration to $\mathcalp^{-1} (V_i)$ and not
only the geometry of the open set $\mathcalp^{-1} (V_i) \subseteq M$ or,
in other words, ``normal forms'' for the singular foliation induced
on $\mathcalp^{-1} (V_i)$ by the restriction of the fibration.
In view of what precedes, the purpose of Theorem~A is to present
explicit models for the tubular neighborhood of a singular fiber
in a fibration $\mathcalp : M \rightarrow S$ whose genus is~$2$. Therefore, except for Section~4, we shall place
ourselves in the following setting: let $D \subset \C$ be the unit disc and suppose that
$\mathcalp$ is a proper holomorphic map from an open complex surface
$M$ to $D$ which satisfies the conditions below:

\noindent 1. $\mathcalp$ defines a regular fibration of
$M \setminus \mathcalp^{-1} (0)$ over $D \setminus \{ 0\}$
with fibers of genus~$2$.

\noindent 2. $\mathcalp^{-1} (0)$ is a connected singular fiber.

To describe the structure of $\mathcalp$ as above, it is convenient to exploit the standard theory of (singular holomorphic)
foliations. In the present case, the foliation in question is given by the fibers of $\mathcalp$. This foliation will be denoted by
$\fol$. In particular the singular points of $\fol$ are exactly the singular points of the singular fiber $\mathcalp^{-1} (0)$.
Our strategy will therefore consist of first describing the structure of $\mathcalp$ on a neighborhood of these singular points.
Then we shall work out the ways in which these local models can be glued together over the singular fiber
$\mathcalp^{-1} (0)$ so as to lead to the full structure of $\mathcalp$ in the sense of Theorem~A.

For the convenience of the reader, we shall first recall some basic notions and index formulas used in the (local) theory
of singular foliations on $(\C^2, 0)$, the reader can consult for example \cite{reference?} for precise definitions, proofs and further information.

In the sequel $\fol$ will stand for a singular holomorphic foliation defined on  a neighborhood of $(0,0) \in \C^2$.
A {\it separatrix}\, for $\fol$ is an irreducible analytic curve passing through the origin and invariant by the foliation.
Suppose that $\fol$ possesses a smooth separatrix $\mathcals$. Modulo changing coordinates, we
can suppose that $\mathcals$ coincides with the axis $\{ y=0 \}$. In these
coordinates, there is a holomorphic vector field $Y$, with isolated singularities and tangent to $\fol$, having the form
$Y = F(x,y) \partial /\partial x + G (x,y) \partial /\partial y$
where, in addition, $G$ is divisible by $y$. The (Camacho-Sad) {\it index}\, of $\mathcals$ with
respect to $\fol$ (at $(0,0) \in \C^2$) is defined by
$$
{\rm Ind}_{(0,0)} (\fol , \mathcals) = {\rm Res}_{x=0}
\frac{\partial}{\partial y} \left( \frac{G}{F} \right ) (x,0) dx \; .
$$
Note, in particular, that ${\rm Ind}_{(0,0)} (\fol, \mathcals) =0$ if $(0,0)$ is a regular point of $\fol$. As it is easy
to see the definition above does not depend on the choices made. The index can be interpreted as measuring
an ``infinitesimal'' self-intersection of $\mathcals$. This interpretation is materialized by the Camacho-Sad formula
as follows. Assume that $\mathcals$ is represented by a global compact Riemann surface $D$ embedded in a complex surface $M$ and invariant by
the foliation $\fol$ (implicitly we now assume that $\fol$ is defined on a neighborhood of $D \subset M$).
Denoting by $p_1 ,\ldots ,p_s$ the singularities of $\fol$ lying
in $D$, we can consider, for each $i =1 , \ldots , s$, the index ${\rm Ind}_{p_i} (\fol, D)$. If $D \cdot D$ stands for
the self-intersection of $D$ then one has
\begin{equation}
\sum_{i=1}^s {\rm Ind}_{p_i} (\fol , D) = D\cdot D \; . \label{camacho}
\end{equation}

The index behaves naturally with respect to blow-ups. Recall
that for a curve $D$ as above, its self-intersection falls by one unity if $D$
is blown-up at a regular point. The index recovers this behavior: if $\tilf$ denotes the blow-up of $\fol$ at the origin,
then the transform $\widetilde{\mathcals}$ of $\mathcals$ is a smooth separatrix for some singularity of $\tilf$. The
index of $\widetilde{\mathcals}$ w.r.t. $\tilf$ equals exactly the index of $\mathcals$ w.r.t. $\fol$ minus $1$. This allows us
to define the index for an irreducible {\it singular separatrix}\, as well. This goes as follows. If $D$ now represents a singular
curve and $p \in D$ is a singular point, then the self-intersection of the blow-up $\widetilde{D}$ of $D$ at $p$ and the self-intersection
of the initial curve $D$ are related by Kodaira's conductor formula. The analogous formula can then be used to define the index
of an irreducible, possibly singular, separatrix for $\fol$. Namely, if $\mathcals$ is a separatrix and $\tilf$ (resp. $\widetilde{\mathcals}$)
stands for the blow-up of $\fol$ (resp. the transform of $\mathcals$), then one has
\begin{equation}
{\rm Ind}_{(0,0)} (\fol , \mathcals) = {\rm Ind}_{q} (\tilf , \widetilde{\mathcals}) + m [\pi^{-1} (0)] \cdot  \widetilde{\mathcals} \; ,
\label{kodairaconductor}
\end{equation}
where $q = \pi^{-1} (0) \cap  \widetilde{\mathcals}$ and where $m$ is the multiplicity of $\pi^{-1} (0)$ as component of
$\pi^{\ast} (\mathcals )$. Finally the intersection product $[\pi^{-1} (0)] \cdot  \widetilde{\mathcals}$ should be regarded
as the usual intersection multiplicity between the curves in question.
Since the transform of $\mathcals$ under a sequence of blow-up transformations will eventually become
smooth, this formula allows us to unequivocally define all these indices. It is rather easy to check that
the resulting index does not depend neither on the number nor on the sequence of blow-up
transformations used to turn $\mathcals$ into a smooth separatrix.

To close this section, let us provide a classification of the singularities of the {\it foliation}\, defined by the fibers of $\mathcalp$ (also said the foliation
associated to $\mathcalp$).
Here it should be noted that, whereas the corresponding singular points coincide with the singular points of the singular
fiber $\mathcalp^{-1} (0)$, the classification mentioned above is not reduced to the description of the singularities of the
corresponding irreducible components of $\mathcalp^{-1} (0)$. To explain this, consider local coordinates $(x,y)$ about a singular
point of $\mathcalp^{-1} (0)$ and let $Y= F(x,y) \partial /\partial x + G (x,y) \partial /\partial y$ be a vector field tangent to fibers of $\mathcalp$
and such that
$F,G$ have only invertible common factors (equivalently the origin is an isolated singularity of $Y$). A {\it normal form}\, for the singularity of the foliation
associated to $\mathcalp$ at the origin
consists of a vector field $Y$ as indicated and of a local holomorphic diffeomorphism taking the (local) orbits of $Y$ to the fibers
of $\mathcalp$ on a neighborhood of the fixed singularity.
This is more restrictive than simply looking at the singularities of the singular fiber viewed as an analytic curve as it will soon be apparent.
First, however, note that the class of singularities of vector fields that appear as singularities of
a singular fibration is very special: it is constituted by those singularities possessing an algebraic holomorphic first integral (unless otherwise
stated, all first integrals mentioned in this work are supposed to be non-constant). By a {\it first integral}\, it is meant a (local) holomorphic function
that is constant over the orbits of the vector field in question (equivalently over the leaves of the associated foliation). Singularities of the
foliation associated to a fibration clearly have a holomorphic first integral which is given by the restriction of the projection $\mathcalp$ itself.
In particular the {\it separatrices}\, of this foliation are contained in the preimage of~$0\in \C$ by the mentioned first integral and thus they
correspond to the (local) irreducible components of the singular fiber $\mathcalp^{-1} (0)$ passing through the singularity in question.

We can now explain the difference between looking at singular points of $\mathcalp^{-1} (0)$ merely as singular points of the singular fiber
or as singular points of the {\it foliation associated}\, to the fibration $\mathcalp$. For example consider the local foliations represented respectively
by the vector fields $Y_1 = x(x-2y) \palix + y(y-2x) \paliy$ and by $Y_2=x(3y-2x) \palix + y(3x-2y) \paliy$. Both singularities appear in the context
of genus~$2$ fibrations (the first one already appear in the case of elliptic fibrations). Their sets of separatrices are isomorphic since they
coincide with the union of the lines $\{ x=0\}, \, \{y=0\}$ and~$\{ x=y\}$. However the foliations associated to these vector fields are not isomorphic in the
sense introduced above since the first has the first integral $xy(x-y)$ and the first integral for the second is $x^2y^2 (x-y)$.

The rest of this section is devoted to providing a complete list of normal forms for the singularities of the {\it foliation associated to $\mathcalp$}.
The first type of singularity to be considered correspond to a {\it nodal singularity}\, of $\mathcalp^{-1} (0)$, i.e. the separatrices of the foliation associated
to $\mathcalp$ in the mentioned singular point consists of two smooth curves with transverse intersection. In suitable coordinates they are respectively
given by $\{ x=0\}, \, \{y=0\}$ so that the corresponding first integral is an invertible multiple of $x^n y^m$, for some $m,n \in \Z_+^{\ast}$. Modulo a new
changing of coordinates, the first integral becomes exactly $x^n y^m$ so that it gives rise to a singularity of the form $mx \palix -ny\paliy$ which will
be called a {\it linear singularity}\, for (the foliation associated
to) $\mathcalp$. Hence every nodal singularity of of $\mathcalp^{-1} (0)$ corresponds to a linear singularity for (the foliation associated to) $\mathcalp$.

In general, given a germ of vector field $Y= F(x,y) \partial /\partial x + G (x,y) \partial /\partial y$ for which $(0,0)$ is an isolated singularity,
the eigenvalues of $Y$ at $(0,0)$ are said to be the eigenvalues of the associated foliation. By construction these are well-defined only up to
a multiplicative constant. In particular, in the case of a linear singularity as above, the eigenvalues of $\mathcalp$ are $m,-n$. Similarly the order of
the foliation at $(0,0)$ (i.e. at the singular point) is the degree of the first non-zero homogeneous component in the Taylor series of $Y$ based at
the $(0,0)$.

The next step is to understand those singularities at which $\mathcalp^{-1} (0)$ locally consists of a number $\alpha \geq 3$ of pairwise transverse smooth curves.
Denote by $\fol$ the local holomorphic foliation induced on a neighborhood of $(0,0) \in \C^2$ by the germ of $\mathcalp$ at this singular point.

\begin{lema}
\label{Section2.11}
Let $\fol$ be as above. Then $\alpha =3$. Besides $\fol$ is conjugate to the foliation given by one of the following first integrals:
$xy(x-y)$, $xy(x-y)^2$, $xy^2(x-y)^3$, $x^2y^2(x-y)$, $x^3y(x-y)$, $x^4y^3(x-y)$, $x^5y^4 (x-y)$ or $x^5y^3 (x-y)^2$.
\end{lema}

\noindent {\it Proof}. Note
that the blow-up $\tilf$ of $\fol$ at the origin will have $\alpha$ singularities over the exceptional divisor $\pi^{-1} (0,0)$
corresponding to the intersections
of $\pi^{-1} (0,0)$ with the transform of each of the (irreducible) separatrices of $\fol$. Besides the first integral of $\fol$ can be written as
$l_1^{k_1}l_2^{k_2}\ldots l_{\alpha}^{k_\alpha} u(x,y)$ with $u(0,0) \neq 0$ and where the $l_i$ are linear forms and the corresponding exponents
$k_i$ are strictly positive integers. Set $k = k_1 + \cdots + k_{\alpha}$ and $P = l_1^{k_1}l_2^{k_2}\ldots l_{\alpha}^{k_\alpha}$. Because the local
leaves of $\fol$ are contained in a global surface of genus~$2$, by considering suitable re-scaling of $\fol$ by means of homotheties of the
form $(x,y) \mapsto (\lambda x ,\lambda y)$ it can be shown that the compactification of the curve $\{P = {\rm cte} \}$ in $\C P (2)$ must be a
curve whose Euler characteristic is not less than~$-2$.

Consider then the curve in $\C^2$ given by $\{ P=1\}$ and identified with its own transform on the blow-up $\widetilde{\C}^2$ of
$\C^2$. The standard projection $\widetilde{\C}^2 \rightarrow \pi^{-1} (0)$, realizes $\{ P=1\}$ as a covering of degree~$k$ of
$\pi^{-1} (0) \setminus \{l_1, \ldots, l_{\alpha}\}$ (where the lines $l_i$ determined by the corresponding linear forms
are identified with their transforms). In particular the Euler characteristic of this affine curve is $k(2-\alpha)$.
Also each singularity $q_i = \pi^{-1} (0) \cap l_i$ of $\tilf$ is linear with eigenvalues $-k_i, k$. Let
$\overline{\mathcal C}$ be the curve obtained by desingularization of the closure of $\{ P=1\}$ in $F_1$ (the blow-up of $\C P (2)$).
Improper points in $\overline{\mathcal C}$ correspond to the improper points of $l_i$, $i=1, \cdots , \alpha$. About these points, there are local coordinates
$(u,v)$ where the curve becomes
$u^{k_i}=v^kF_i(u,v)$
with $F_i (0,0) \neq 0$. It follows that the number of branches at this point is the greatest common divisor (g.c.d.) between $k_i$ and $k$. Hence the Euler
characteristic $\chi(\overline{{\mathcal C}})$ of $\overline{\mathcal C}$ is obtained by adding one unit for branch at infinity to the Euler characteristic of $\{ P=1\}$.
Therefore we have:
$$
\chi(\overline{{\mathcal C}})= k(2-\alpha) + \sum_{i=1}^{\alpha} {\rm g.c.d.}\, (k_i , k).
$$
Thus
$$
k(2-\alpha)+\alpha \leq 
\chi(\overline{{\mathcal C}})
\leq
k(2-\alpha) + \sum_{i=1}^{\alpha}{k_i} =k(3-\alpha).
$$
However it was seen that the Euler charcteristic of $\overline{\mathcal C}$ must take on the values $2, 0, -2$. If it were equal to~$2$, it would follow that
$\alpha =2$ what is impossible since $\alpha \geq 3$ (clearly the case $\alpha =2$ corresponds to a linear singularity). Next suppose that
$\chi(\overline{{\mathcal C}})=0$. In this case the above estimates show that $\alpha =3$ and that ${\rm g.c.d.}\, (k_i , k) =k_i$ for $i=1,2,3$. It is
immediately to check that the only possibilities are
$$
(k_1,k_2,k_3) = (1,1,1), (1,1,2), (1,2,3).
$$
Finally suppose that $\chi (\overline{\mathcal C})=-2$. Again it follows that $\alpha =3$ since $k > 2$. Thus $k = 2 + \sum_{i=1}^{3} {\rm g.c.d.}\, (k_i , k)$
with $k_1 + k_2 + k_3 =k$ and ${\rm g.c.d.}\, (k_1, k_2, k_3) =1$. The new non-trivial solutions occur for $k=5,\, 8, \, 10$ and they correspond to
$$
(k_1,k_2,k_3) = (2,2,1), (3,1,1), (4,3,1), (5,4,1), (5,3,2)\, .
$$
The preceding provides us with normal forms for $P$. To show that $\fol$ is actually conjugate to the foliation $\fol_P$ induced by the corresponding $P$
is now very standard and goes as follows. Since $\fol$ (resp. $\fol_P$) possesses only $3$ separatrices $d^1, d^2, d^3$ (resp. $d_P^1, d_P^2, d_P^3$)
that, in addition, are smooth a pairwise transverse, it is well-known the existence of a local holomorphic fibration ${\rm Fib}$ (resp. ${\rm Fib}_P$) defined on a
neighborhood of the exceptional divisor $\pi^{-1} (0,0)$ and satisfying the following conditions:
\begin{itemize}
\item The transforms of $d^1, d^2, d^3$ (resp. $d_P^1, d_P^2, d_P^3$) are contained in fibers of ${\rm Fib}$ (resp. ${\rm Fib}_P$).

\item Away from the above mentioned fibers, ${\rm Fib}$ (resp. ${\rm Fib}_P$) is transverse to the blow-up $\tilf$ (resp. $\widetilde{\fol}_P$)
of $\fol$ (resp. $\fol_P$).
\end{itemize}
On the other hand, the lines $d^1, d^2, d^3$ and $d_P^1, d_P^2, d_P^3$ determine the same set of singular points $q_1, q_2, q_3$ in
$\pi^{-1} (0,0)$ for both foliations $\fol$ and $\fol_P$. In other words, $\pi^{-1} (0,0) \setminus \{q_1, q_2, q_3 \}$ is a common leaf for both
foliations $\fol, \, \fol_P$. The holonomy representations of this leaf with respect to $\fol$ and with respect to $\fol_P$ are holomorphically
conjugate. In fact, they are both holomorphically conjugate to the same finite group of rotations, cf. Lemma~\ref{lema4.1} for a general argument.
Now a holomorphic conjugacy between $\tilf, \, \widetilde{\fol}_P$ can be constructed on a neighborhood of $\pi^{-1} (0,0)$ by the standard
method of ``lifting paths'' thanks to the existence of the fibrations ${\rm Fib}, \, {\rm Fib}_P$ and to the fact the common leaf $\pi^{-1} (0,0) \setminus \{q_1, q_2, q_3 \}$
induce conjugate holonomy representations. This completes the proof of the lemma.\qed

Again let $\fol$ denote the local holomorphic foliation induced on a neighborhood of $(0,0) \in \C^2$ by the germ of (the foliation associated
to) $\mathcalp$ at a singular point. The next proposition provides normal forms for all these local foliations. To establish this proposition it is convenient to
refer to linear singularities as singularities of {\it generation zero}\, and to the (non-linear) singularities appearing in Lemma~\ref{Section2.11}
as singularities of {\it generation one}. In the proof of the proposition below the reader is assumed to be familiar with Seidenberg theorem in \cite{sei}.

\begin{prop}
\label{Section2.22}
If $\fol$ is not conjugate to any of the models presented above, it must be conjugate to the foliation given by one of the following first integrals:
$y(x^2-y)$, $y^2(x^2-y)$, $y^3 (x^2-y)$, $y^4 (x^2-y)$, $y^3 (x^2-y)^2$, $x^2 - y^3$, $y^2 (x^2-y^3)$ and $x^2 -y^5$.

\end{prop}

\noindent {\it Proof}. We blow-up $\fol$ and look at the singularities of $\tilf$ over $\pi^{-1} (0)$. If all these singularities are linear, then we have
one of the possibilities indicated above. Thus let us now assume that these singularities are either linear or such that all its (irreducible) separatrices
are smooth and pairwise transverse (in other words, these singularities are of generation zero or one).
Naturally it is also assumed that at least one of these singularities, say $q_1$, is of generation one.
A simple application of the index formula~(\ref{camacho}) shows that $q_1$ must be $\pi^{-1} (0)$ the unique singularity of $\tilf$ in
$\pi^{-1} (0)$ (and its separatrix defined by $\pi^{-1} (0)$ must have index equal to~$-1$). In fact, the indices associated to separatrices
in all the singularities above is strictly negative. Besides these indices are always less than or equal to~$-1$ for generation one singularities. Collapsing
$\pi^{-1} (0)$ leads to the second generation of singularities, namely those of the form: $y(x^2-y)$, $y^2(x^2-y)$, $y^3 (x^2-y)$, $y^4 (x^2-y)$ and
$y^3 (x^2-y)^2$.

We continue by recurrence, we consider now that $\pi^{-1} (0,0)$ contains singularities of generations zero, one and two with at least one singularity $q_1$
of generation two. The previous remark concerning the fact that all separatrices have indices strictly negative is still verified by singularities of
generation two. In particular it follows that $\pi^{-1} (0,0)$ cannot contain singularities of generation one. Now a direct inspection on the indices of separatrices
for singularities of generation two yields the following:
\begin{itemize}

\item the indices of the separatrices of $y(x^2-y)$ and $y^3 (x^2-y)^2$ take only on values strictly smaller
than~$-1$. Therefore these singularities cannot appear on $\pi^{-1} (0,0)$. Therefore these singularities are terminal in the Seidenberg procedure.

\item $y^2(x^2-y)$ has a unique separatrix of index~$-1$. Its collapsing leads to the model $x^2 - y^3$ with a single separatrix which happens to be
singular. In particular this singularity cannot be combined with any other over a rational curve of self-intersection~$-1$, therefore it is also terminal
in the Seidenberg procedure.

\item The separatrix $\{ y=0\}$ of $y^3 (x^2-y)$ (resp. $y^4 (x^2-y)$) has index~$-2/3$ (resp. $-1/2$) while the other separatrix has index
$-6$ (resp. $-8$).
\end{itemize}
Besides the index formula, there are some elementary relations that must be verified by germs of singular vector fields appearing in the exceptional
divisor obtained by blowing-up another (holomorphic) vector field. In view of these relations, it is easy to see that
the singularity $y^3 (x^2-y)$ cannot be combined with another singularity as above over a rational curve of self-intersection~$-1$. To check this claim
it is enough to use the fact that the sum of the ``asymptotic orders'' of singularities lying over a rational curve equals~$2$ (see \cite{re1} for a very detailed discussion).
This singularity is therefore terminal as well. Finally, with the same type of arguments, it follows that the singularity $y^4 (x^2-y)$ can be combined
(over a rational curve of self-intersection~$-1$) with a linear singularity of  eigenvalues~$-1,2$. This combination yields the model associated to the first
integral $y^2 (x^2 - y^3)$ which is the only singularity of generation three. In turn this singularity has a singular separatrix (that therefore cannot be contained in
a exceptional divisor) and another smooth separatrix of self-intersection~$-1$. Thus it must be the unique singularity over a
rational curve of self-intersection~$-1$. The collapsing of this rational curve leads to the generation four singularity $x^2 -y^5$ which has a unique separatrix
that happens to be singular. In particular this singularity must be terminal. The proposition is proved.\qed

To close this section let us indicate how the singularities identified above are related to some special isotrivial fibrations of genus~$2$.
We shall do it only in the case of the singularity $x^5 y^4(x-y)$ since the remaining cases are analogous. The vector field (foliation) associated to this
first integral is $x(4x-5y) \palix - y(6x-5y) \paliy$. In particular it defines a foliation on all of $\C^2$ and, in fact, on $\C P (2)$. This foliation will be
denoted by $\fol$ and the reader will notice that the ``line at infinity'' $\Delta$ is invariant by $\fol$. The lines
$\{x = 0\}$, $\{y = 0\}$ and $\{x = y\}$ are invariant by $\fol$ and meet $\Delta$ at points denoted respectively by
$P_1, \, P_2$ and $P_3$. Thus, on $\C P(2)$, the corresponding foliation (pencil) leaves $4$ projective lines invariants and possesses $4$ singularities,
namely the origin $(0,0) \in \C^2 \subset \C P(2)$ and the points $P_1, \, P_2$ and $P_3$.

\begin{figure}[hbtp]
\centering
\psfrag{Y}{$Y$}
\psfrag{X}{$X$}
\psfrag{Z}{$Z$}
\psfrag{De}{$\Delta$}
\psfrag{P1}{$P_1$}
\psfrag{P2}{$P_2$}
\psfrag{P3}{$P_3$}
\includegraphics[scale=0.3]{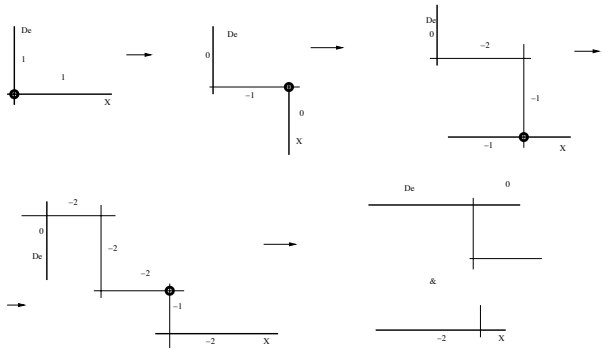}
\caption{Blow up of first singularity}
\label{N6P3blow}
\end{figure}

On a neighborhood of $P_3$, we write $x = \frac{1}{u}$, $y = \frac{v}{u}$. The first integral
becomes $F(u, v) = \frac{v^4(1-v)}{u^{10}}$ so that  our singularity looks like $\frac{v^4}{u^{10}}$
near $u = v= 0$. Therefore this first integral has a ``point of indetermination'' at $P_3$ which, in turn, can be
eliminated by an appropriate sequence of blowing-ups. This sequence  is depicted below (Figure~\ref{N6P3blow}),
each of the numbers on the components indicates its self-intersection.

\begin{figure}[hbtp]
\centering
\psfrag{Y}{$Y$}
\psfrag{X}{$X$}
\psfrag{Z}{$Z$}
\psfrag{De}{$\Delta$}
\psfrag{B}{$\text{and}$}

\includegraphics[scale=0.3]{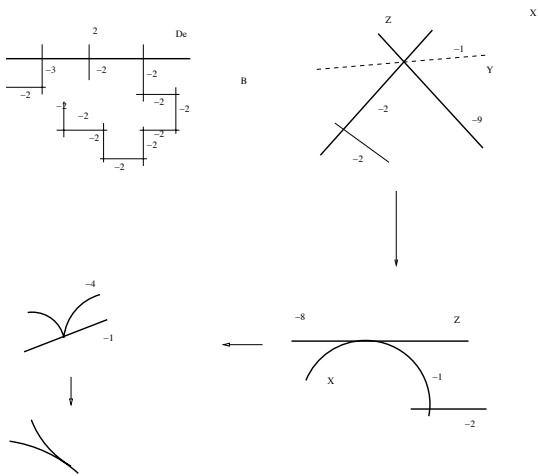}
\caption{Type $20$ and a pinched torus}
\label{N6}
\end{figure}

\noindent An analogous analysis can be carried out at singularities $P_2$ and $P_1$.
Gathering all the information about these three singularities, we obtain the following dual fibers, 
which are classified as Type $20$ (since it has a component of self-intersection $3$) 
and a pinched torus in \cite{ogg} (note that  in Figure~\ref{N6} the components with 
self-intersection $-1$ were collapsed).

\section{Proof of Theorem~A}

The goal of this section is to prove Theorem~A which can be thought of as a local uniformization theorem for
fibrations of genus~$2$. The strategy of the proof relies, on one hand, on the study of the singularities of the fibration (viewed as foliation) for
each model in Ogg's list and, on the other hand, on the analysis of the ``total holonomy
group of the singular fiber''.
In particular in the case of nodal singularities, it will be shown that the corresponding eigenvalues are uniquely determined by the combinatorial data
of the singular fiber itself.

Concerning the index of a separatrix with respect to a singular foliation $\fol$, as defined in Section~$2$, let us consider the special case
in which the foliation $\fol$ possesses non-zero eigenvalues at the singularity. More precisely, suppose that at the origin $\fol$ has eigenvalues
$\lambda_1 , \, \lambda_2 \neq 0$ with
quotient belonging to $\C \setminus \R_+$. By definition this means that $\fol$ can be represented by a local holomorphic
vector field whose linear part at the origin has $\lambda_1 , \, \lambda_2$ as eigenvalues. It is well-known that, in this case,
there are local coordinates $(x,y)$ where $\fol$ is represented by a vector field of the form
$$
\lambda_1 x (1 + {\rm h.o.t.})\frac{\partial}{\partial x} 
+ 
\lambda_2 y (1 + {\rm h.o.t.})\frac{\partial}{\partial y} \, .
$$
In particular $\fol$ possesses exactly two separatrizes given in the above coordinates by the axes
$\{ x=0 \}$ and $\{ y=0 \}$. A direct inspection shows that
$$
{\rm Ind}_{(0,0)} (\fol , \{ y = 0 \}) = \frac{\lambda_2}{\lambda_1}
\hspace{1cm} {\rm and } \hspace{1cm}
{\rm Ind}_{(0,0)} (\fol , \{ x = 0 \}) = \frac{\lambda_1}{\lambda_2}.
$$
Hence
\begin{equation}
{\rm Ind}_{(0,0)} (\fol , \{ y = 0 \}) = \frac{1}{{\rm Ind}_{(0,0)} (\fol , \{ x = 0 \})} \, . \label{inversion}
\end{equation}

It the context of fibrations it is particularly easy to detect when a singular point of a fibration $\fol$ possesses non-zero eigenvalues: this happens
if and only if the point in question is a nodal singularity of the singular fiber (cf. Section~2). It the follows that at these points the singularity of
$\mathcalp$ (viewed as singularity of a foliation)
is totally determined by its eigenvalues. Based on this remark, we are going to prove the following stronger
statement. Recall that the Dynkin diagram of a singular fiber is supposed to include local multiplicites (or their
self-intersection numbers cf. Section~2).

\begin{prop}
\label{localmodels}
Let $\mathcalp : M \rightarrow D\subset \C$ denote a fibration of genus~$2$ having $\mathcalp^{-1} (0)$ as singular
fiber. Then the Dynkin diagram of $\mathcalp^{-1} (0)$ analytically
determines the structure of the fibration $\mathcalp$ on a neighborhood of every singular point.
\end{prop}

The proposition above can also be stated as follows.
Suppose that $\mathcalp_1 : M \rightarrow D\subset \C$ (resp. $\mathcalp_2 : M \rightarrow D\subset \C$) are fibrations
of genus~$2$ whose singular fibers $\mathcalp_1^{-1} (0), \; \mathcalp_2^{-1} (0)$ have isomorphic Dynkin diagrams.
Then, if $p_1 \in \mathcalp_1^{-1} (0)$, $p_2 \in \mathcalp_2^{-1} (0)$ are corresponding singularities in this Dynkin diagram,
there is a local holomorphic diffeomorphism from $p_1$ to $p_2$ sending fibers of $\mathcalp_1$ to fibers
of $\mathcalp_2$.

\bigskip

\noindent {\it Proof of Propostion~\ref{localmodels}}. Let us start with the case of fibers having only nodal singularities.
As already seen, on a neighborhood of a nodal singularity, the fibration $\mathcalp$ is given by a linear vector field
$mx \palix - ny \paliy$, $m,n \in \N^{\ast}$. In particular the indices of each axis $\{ x=0 \}$, $\{ y=0 \}$ with respect to
$\mathcalp$ verify Equation~\ref{inversion}. All we have to do is to check that the eigenvalues of each singularity  of
$\mathcalp^{-1}(0)$ is then determined by its position in the corresponding Dynkin diagram.

\begin{figure}[hbtp]
\centering
\psfrag{p1}{$p_1$}
\psfrag{p2}{$p_2$}
\psfrag{p3}{$p_3$}
\psfrag{p4}{$p_4$}
\psfrag{p5}{$p_5$}
\psfrag{p6}{$p_6$}
\psfrag{S1}{$S_1$}
\psfrag{S2}{$S_2$}
\psfrag{S3}{$S_3$}
\psfrag{S4}{$S_4$}
\psfrag{S5}{$S_5$}
\psfrag{S6}{$S_6$}
\psfrag{S7}{$S_7$}

\includegraphics[scale=0.3]{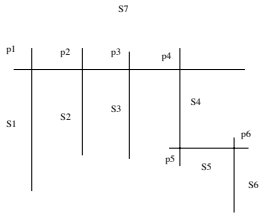}
\caption{Description of the eigenvalues for Model $4$ }
\label{Evals_Model4}
\end{figure}

The easiest way to deal with this question is to check each model with linear singularities individually.
To avoid unnecessary repetition, we shall explain explicitly some examples
that contain the general procedure.
First let us consider the model of singular fiber noted {\bf Model $4$} in \cite{ogg}, cf. Figure~\ref{Evals_Model4}.

\noindent We begin by the components that only intersect the rest of the fiber once.
\begin{equation}
{\rm Ind}_{p_1} (\fol , S_1) = S_1 \cdot S_1 = -4 =  \frac{\lambda_{2,1}}{\lambda_{1,1}}
\hspace{0.3cm} {\rm and} \hspace{0.3cm}
{\rm Ind}_{p_i} (\fol , S_i) = S_i \cdot S_i = -2 =  \frac{\lambda_{2,i}}{\lambda_{1,i}},
\end{equation}
for $i = 2, 3, 6$. Next note that, by changing variables (or re-scaling coordinates), we can always take one of the eigenvalues of this type of components (say, for example, 
$\lambda_{1,i}$, $i = 1, 2, 3, 6$) to be equal to $1$. 
Hence, we determine that $\lambda_{2,1} = -4$, and $\lambda_{2, i} = -2$, for $i = 2, 3, 6$. Also, 
\begin{eqnarray}
{\rm Ind}_{p_5} (\fol , S_5) + {\rm Ind}_{p_6} (\fol , S_5) & = & 
\frac{1}{\lambda_{2,6}} + \frac{\lambda_{2,5}}{\lambda_{1,5}}  = -2
\\
\sum_{j=1}^4  {\rm Ind}_{p_j} (\fol , S_7) & = & \sum_{j=1}^4  \frac{\lambda_{2,j}}{\lambda_{1,j}} = -2.
\end{eqnarray}
The first equation determines ${\rm Ind}_{p_5} (\fol , S_5) = -\frac{3}{2}$, 
and the last determines ${\rm Ind}_{p_4} (\fol , S_7) = -\frac{3}{4}$.

Note that we would still have an extra relation given by the separatrix $S_4$, which is redundanct.
This phenomenon will occur in all examples where the fiber model corresponds to a graph with trivial
fundamental group. It is due to the trivial fact that, in such graphs, the number of vertices is
strictly smaller than the number of edges.

\begin{figure}[hbtp]
\centering
\psfrag{p1}{$p_1$}
\psfrag{p2}{$p_2$}
\psfrag{p3}{$p_3$}
\psfrag{p_4}{$p_4$}
\psfrag{p5}{$p_5$}
\psfrag{p6}{$p_6$}
\psfrag{S1}{$S_1$}
\psfrag{S2}{$S_2$}
\psfrag{S3}{$S_3$}
\psfrag{S4}{$S_4$}
\psfrag{S5}{$S_5$}
\psfrag{S6}{$S_6$}
\psfrag{S7}{$S_7$}

\includegraphics[scale=0.4]{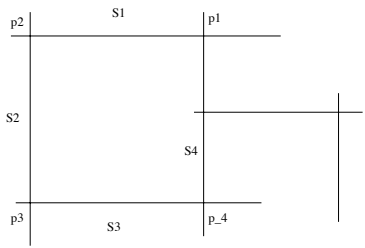}
\caption{Description of the eigenvalues for Model $10$ }
\label{Evals_Model10}
\end{figure}

Still considering only linear singularities, there is another possible case in which some of the components
of the singular fiber form a loop. As a prototype for these cases, let us consider
Model $10$ in \cite{ogg}, cf. Figure~\ref{Evals_Model10}.
Following the same reasoning as above, and writing ${\rm Ind}_{p_i} (\fol , S_j) = I_{i,j}$,
we see that
$$
I_{1,1} + I_{4,1} = -\frac{4}{3}; \hspace{0.5cm} 
I_{1,2} + I_{2,2} = -2; \hspace{0.5cm}
I_{2,3} + I_{3,3} = -4; \hspace{0.5cm} {\rm and} \hspace{0.5cm}
I_{3,4} + I_{4,4} = -2.
$$
This $4 \times 8$ system has a rank $4$ space of solutions, which will give us uniqueness 
when we impose the conditions that, if two edges $S_i$ and $S_j$ intersect, then for any 
$p_q \in S_i \cap S_j$, $I_{q,i} I_{q,j} = 1$.

Once again, this example reflects the general idea: whenever there is a loop in the graph representing
the fiber, the number $N$ of edges will necessarily match the number of vertices.
Camacho-Sad index formula will provide a $N \times 2N$ linear system of rank $N$ which combined to
the observation on the previous paragraph implies the uniqueness of the indices of the singularities
for each model.
In general the study of the singularities for each of the models in \cite{ogg} that have 
linear singularities will coincide with one of the two cases depicted above, and
therefore is going to be omitted.

It remains to discuss the cases corresponding to the singular fibers exhibiting non-linear singularities. As previously
seen, non-linear singularities are detected as the singular points of $\mathcalp^{-1} (0)$ that are not of nodal type. Normal forms
for (the foliation associated to) $\mathcalp$ on a neighborhood of a non-linear singularity were presented in Section~2.2.
It can directly be checked that each possible normal form is totally identified
by the nature of the singularity of $\mathcalp^{-1} (0)$ (as a singularity of a curve) along with the corresponding indices.
For example suppose that $\mathcalp^{-1} (0)$ is locally given by three smooth curves pairwise transverse and with indices
$-3/2, -3/2, -4$. Though all the models listed in Proposition~\ref{Section2.11} are such that their corresponding separatrices
consist of three smooth curves pairwise transverse, there is only one whose separatrices match the mentioned indices,
namely the vector field corresponding the the first integral $x^2 y^2 (x-y)$. In other words, whenever we reach a singularity
of $\mathcalp^{-1} (0)$ satisfying the above conditions, the local structure of $\mathcalp$ around this singularity is determined by
the first integral $x^2y^2 (x-y)$. A similar statement hold for the models appearing in Proposition~\ref{Section2.22}.

Summarizing to deal with the cases where there are non-linear singularities we proceed as follows. By considering
the linear singularities appearing in the model in question, we repeat the preceding
analysis to conclude that the indices of the separatrizes of the {\it non-linear}\, non-linear singularity are uniquely
determined. By the above observation, these indices together with the form of the singularity formed by the
separatrizes themselves characterize unequivocally the local normal form of $\mathcalp$. This concludes the proof of
the proposition.\qed

What precedes provides a complete description of the structure of the fibration $\mathcalp$ on a neighborhood of a
singular point of $\mathcalp^{-1} (0)$. Thus, by now, we have understood the ``lomathcal Pieces'' that can be used to build
$\mathcalp$. Our next task will consist of working out the possible assembling of these pieces so as to arrive to the proof of
Theorem~A.
The key notion that will lead us to the proof of this theorem is the holonomy associated to $\mathcalp^{-1} (0)$. To explain
this notion, let us fix an irreducible component $C$ of $\mathcalp^{-1} (0)$. Then $L = C \setminus {\rm Sing}\, (\mathcalp^{-1} (0))$
can be viewed as a regular leaf of the foliation induced $\mathcalp$. Let then $\Sigma$ denote a local transverse
section to $L$ at a point $x \in L$. In particular we note that every fiber of $\mathcalp$ can cut $\Sigma$ at a uniformly bounded
number of points. Identifying $x, \, \Sigma$ with a neighborhood of $0 \in \C$, the holonomy of $L$ provides us
a representation
$$
\rho_C : \, \pi_1 (L) \longrightarrow {\rm Diff}\, (\C ,0)
$$
where $\pi_1 (L)$ stands for the fundamental group of $L$.

\begin{lema}
\label{lema4.1}
The image $\rho_C (\pi_1 (L)) \subset {\rm Diff}\, (\C ,0)$ is a finite abelian group. In suitable coordinates
it is generated by a rational rotation.
\end{lema}

\noindent {\it Proof}.  To show that $\rho_C (\pi_1 (L))$ is abelian, let us suppose for a contradiction that it is not
the case. Hence there are elements $h_1, h_2 \in \rho_C (\pi_1 (L))$ such that $h = h_1 \circ h_2 \circ h_1^{-1}
\circ h_2^{-1}$ is not reduced to the identity. However the resulting local diffeomorphisms $h \neq {\rm Id}$
obviously satisfies $h'(0) =1$. The local dynamics of this type of diffeomorphism is known as the ``flower'' and it
possesses infinite orbits accumulating at $0 \in \C$. This means that nearby fibers of $\mathcalp$ would accumulate
on $\mathcalp^{-1} (0)$ what is impossible.

Let now $h\neq {\rm Id}$ be an arbitrary element of $\rho_C (\pi_1 (L))$. Clearly we must have $\vert h'(0) \vert
=1$ since otherwise $h$ or $h^{-1}$ would have orbits accumulating at $0 \in \C$ what is impossible. Hence we
can set $h(z) = e^{2\pi i \theta} z + \cdots$. If $\theta$ is not rational then we would still have fiber of $\mathcalp$ intersecting
$\Sigma$ a number arbitrarily large of times. As already seen, this is again impossible. Thus we finally conclude that
$h(z) = e^{2\pi i p/q} z + \cdots$. In particular $h^q (z) = z + \cdots$. In view of the preceding argument, it follows
that $h^q = {\rm Id}$.

The preceding shows that $\rho_C (\pi_1 (L))$ is a finitely generated abelian group all of whose elements are
of finite order. It follows then that $\rho_C (\pi_1 (L))$ is finite and hence linearizable thanks to the standard B\"ochner theorem.
We then conclude that $\rho_C (\pi_1 (L))$ is conjugate to a (finite) group of rational rotations which, in turn, must be cyclic.
The lemma is proved.\qed

Consider local genus~$2$ fibrations $\mathcalp_1 : M_1 \rightarrow D$, $\mathcalp_2 : M_2 \rightarrow D$ as before. Now that non-linear singularities of $\mathcalp$ are understood
we can drop the assumption of having minimal singular fibers i.e. $\mathcalp_1^{-1} (0), \, \mathcalp_2^{-1} (0)$ are allowed to contain rational curves with
self-intersection equals to~$-1$. The advantage of doing so is that all the singularities of $\mathcalp_1, \, \mathcalp_2$ become linear and we still keep
an obvious correspondence between the combinatorial data of $\mathcalp_1^{-1} (0), \, \mathcalp_2^{-1} (0)$. Naturally without the preceding proposition we would not
be able to assert that a correspondence between ``minimal models'' of $\mathcalp_1^{-1} (0), \, \mathcalp_2^{-1} (0)$ implies a similar correspondence for the non-minimal
models obtained by proceeding further blow-ups aiming at linear singularities.

Finally a simple but central observation to prove Theorem~A is as follows:

\begin{lema}
\label{lema4.22}
Let $C$ be an irreducible component of $\mathcalp_1^{-1} (0)$ and denote by $C_1, \ldots , C_r$ the other irreducible
components of $\mathcalp_1^{-1} (0)$ intersecting $C$. Then there exists a neighborhood $U$ of $C$ equipped with a
$C^{\infty}$-fibration $\xi_{1,C} : U \rightarrow C$ such that the following holds:
\begin{itemize}
\item The intersections $C_i \cap U$ are connected and contained in fibers of $\xi_{1,C}$, $i=1, \ldots, r$.

\item Away from $C_i \cap U$, $i=1, \ldots, r$, the fibers of $\xi_{1,C}$ are transverse to the fibers of $\mathcalp$.
\end{itemize}
\end{lema}

\noindent {\it Proof}. Recall the singularities of $\mathcalp_1^{-1} (0)$ are now linear. These singularities are in natural correspondence
with the intersection points $p_i = C \cap C_i$, $i=1, \ldots, r$. The existence of a local fibration with the desired properties on a neighborhood
of $p_i$ is therefore obvious since $\mathcalp$ is locally conjugate to a model of the form $mx \palix - ny\paliy$, $m,n \in \Z_+$.

On the other hand on a compact part $K$ of $C \setminus \{p_1, \ldots , p_r \}$ the statement is an immediate consequence of the standard
$C^{\infty}$ theorem of tubular neighborhood. Now by using bump functions it is easy to check that the latter fibration can be glued together with the above
mentioned local fibrations on neighborhoods of the singular points $p_i$, $i=1, \ldots, r$. This proves the lemma.\qed

We are now ready to start the approach to the proof of Theorem~A.

\begin{prop}
\label{quasithm}
Consider fibrations $\mathcalp_1, \, \mathcalp_2$ as above and suppose that all the irreducible components of $\mathcalp_1^{-1} (0), \, \mathcalp_2^{-1} (0)$
are rational curves. Suppose also that the Dynkin diagrams associated to these singular fibers contain no loop. Then there exists
a transversely holomorphic $C^{\infty}$-conjugacy between $\mathcalp_1, \, \mathcalp_2$.
\end{prop}

\noindent {\it Proof}. Let $C^1$ be an irreducible component of $\mathcalp_1^{-1} (0)$ and denote by $C^2$ the corresponding
component of $\mathcalp_2^{-1} (0)$. Denote by $p_1^1, \ldots , p_r^1$ the singularities of $\mathcalp_1$ lying in $C^1$. Each $p_i^1$
is a linear singularity whose local holonomy is conjugate to a rational rotation. Fixed a base point $x_0 \in C^1 \setminus
\{ p_1^1, \ldots , p_r^1 \}$, the fundamental group of $C^1 \setminus \{ p_1^1, \ldots , p_r^1 \}$ is generated by loops $c_1,
\ldots , c_k$ based at $x_0$ and each of them encircling a single singularity $p_1^1, \ldots , p_r^1$. In particular $c_i$ is
freely homotopic to a small circle about $p_i$ characterizing the local holonomy of this singularity. The only non-trivial
relation verified by the generators $c_1, \ldots , c_k$ is $c_1 \cdot c_2 \cdots c_k ={\rm id}$. Fixed a local transverse
section $\Sigma$ at $x_0$, the holonomy of the regular leaf of $\mathcalp_1$ given precisely by $C^1 \setminus \{ p_1^1, \ldots , p_r^1 \}$
can be identified to a homomorphism
$$
\rho : \pi_1 (C^1 \setminus \{ p_1^1, \ldots , p_r^1 \} ) \longrightarrow {\rm Diff}\, (\C, 0) \, .
$$
Recalling that $\rho (c_i) = h_i$ is conjugate to a rational rotation, we can denote its order by $n_i \in \N^{\ast}$. Thus
the image $\rho (\pi_1 (C^1 \setminus \{ p_1^1, \ldots , p_r^1 \} )) \subset {\rm Diff}\, (\C, 0)$ is generated by local
diffeomorphisms $h_1, \ldots , h_r$ satisfying the following relations:
\begin{equation}
h_1^{n_1} = \cdots = h_r^{n_r} = h_1 \circ  \cdots \circ h_r   = {\rm id} \; .
\end{equation}
The argument used in the proof of Lemma~\ref{lema4.1} still implies that $\rho (\pi_1 (C^1 \setminus \{ p_1^1, \ldots , p_r^1 \} ))$
is abelian. Thus it is also finite and generated by a single rational rotation whose order is totally determined by the orders
$n_1, n_2, \ldots ,n_r$. In turn this means that $\rho (\pi^1 (C^1 \setminus \{ p_1^1, \ldots , p_r^1 \} ))$ is generated by a
rational rotation whose order is determined explicitly by the eigenvalues associated to the singularities
$p_1^1, \ldots , p_r^1$.

Let us now consider the corresponding components
$C^1,C^2$ along with the singularities $\{ p_1^1, \ldots , p_r^1 \} \in C^1$ (resp. $\{ p_1^2, \ldots , p_r^2 \}
\in C^2$). According to Proposition~\ref{localmodels}, for every $i=1, \ldots ,r$ the eigenvalues of
$\mathcalp_1$ at $p_i^1$ coincide with those of $\mathcalp_2$ at $p_i^2$. Then what precedes ensures us that the
holonomy of the leaf $C^1 \setminus \{ p_1^1, \ldots , p_r^1 \}$ w.r.t. $\mathcalp_1$ is analytically conjugate to the holonomy
of $C^2 \setminus \{ p_1^2, \ldots , p_r^2 \}$ w.r.t. $\mathcalp_1$ since they are both conjugate to the group generated by the same rational
rotation. With this information in hand, we can proceed to construct a conjugacy between $\mathcalp_1$ and
$\mathcalp_2$ on a neighborhood of $C^1, C^2$ as follows. For each $i=1, \ldots ,r$ and $j=1,2$, let $W_i^j$ denote a
small neighborhood of $p_i^j$. In particular $C^1 \setminus \bigcup_{i=1}^r W_i^1$ (resp.
$C^2 \setminus \bigcup_{i=1}^r W_i^2$) is a compact part of the leaf $C^1 \setminus \{ p_1^1, \ldots , p_r^1 \}$
(resp. $C^2 \setminus \{ p_1^2, \ldots , p_r^2 \}$). The holomorphic conjugacy between the corresponding holonomy groups can then
be extended to a (transversely holomorphic) $C^{\infty}$-conjugacy between $\mathcalp_1, \mathcalp_2$ on neighborhoods
of $C^1 \setminus \bigcup_{i=1}^r W_i^1$, $C^2 \setminus \{ p_1^2, \ldots , p_r^2 \}$. Naturally this extension is by the standard method
of ``lifting paths'' which, in turn, depends on the transverse fibration constructed in Lemma~\ref{lema4.22}. Because this fibration is of
class $C^{\infty}$ this guarantees the $C^{\infty}$-character of the resulting conjugacy. Finally the fact that
the singularities $p_i^j$ are linear and have a ``saddle-like behavior'' (a consequence of the sign of the eigenvalues in question) makes it
easy to check that this conjugacy can be extended to the neighborhoods $W_i^j$.

Next we need to show that the conjugacy
constructed above can be extended from the component $C^1 \in \mathcalp_1^{-1} (0)$ (resp. $C^2 \in \mathcalp_2^{-1} (0)$)
to a subsequent component of $\mathcalp_1^{-1} (0)$ (resp. $\mathcalp_2^{-1} (0)$). This goes as follows. First note that
the holonomy group of $C^1$ (resp. $C^2$) w.r.t. $\mathcalp_1$ (resp. $\mathcalp_2$) as defined above may differ from a similar
notion of holonomy that takes into account the entire singular fiber $\mathcalp_1^{-1} (0)$. To explain this difference, the
holonomy representation considered above for the component $C^1$ will be referred to as the {\it holonomy generated at $C^1$}.
Next consider another rational curve $C_2^1$
contained in $\mathcalp_1^{-1} (0)$ and intersecting $C^1$ at a linear singularity $p$. Denote by $m,-n$ the eigenvalues
of $\mathcalp_1$ at $p$ so that the local holonomy $h_1^1$ of $C^1$ (resp. $h_2^1$ of $C_2^1$) is conjugate to a rotation of order $m$
(resp. $n$). The singularity $p$ can allow part of the holonomy of $C^1$ to be ``transmitted'' to $C_2^1$ (and vice-versa).
Indeed, the effect of the singularity $p$ is to change a rotation of order $m$ into a rotation of order $n$. Thus, if the
holonomy group generated at $C^1$ has order exactly $m$ then there is no transmission from $C^1$ to
$C_2^1$. However if the group in question has order strictly larger than $m$, then those transformations that are not
in the subgroup generated by $h_1^1$ will induce non-trivial transformations on $C_2^1$. These new transformations
need not be contained in the {\it holonomy group generated at $C_2^1$}. Finally, the group of transformations induced at
an irreducible component $D_i^1$ by all these transformations will be called {\it the total holonomy of $C_i^1$ w.r.t.
$\mathcalp_1$}. The argument given above actually shows that every two corresponding irreducible
components $C^1 \in \mathcalp_1^{-1} (0)$ and $C^2 \in \mathcalp_2^{-1} (0)$ have conjugate total holonomy. It is now easy to see
that the above constructed conjugacy on a neighborhood of $C^1$ can be extended to the other components of this singular
fiber. In turn these extensions are well-defined since the Dynkin
diagram of the singular fibers contains no loop. The proposition is proved.\qed

\vspace{0.2cm}

\noindent {\it Proof of Theorem~A}. Keeping the above setting, let us first consider the case where the irreducible components
of $\mathcalp_1^{-1} (0), \, \mathcalp_2^{-1} (0)$ are all rational curves. The corresponding Dynkin diagram is however allowed to possess
loops. In view of Section~2.1, these loops are entirely constituted by rational curves. Besides the singular fibers may contain at most
two of these loops. Fix a loop contained in $\mathcalp_1^{-1} (0)$ (resp. $\mathcalp_2^{-1} (0)$) and note that this loop contains a ``core curve'' $c$ that
is homotopically non-trivial. Therefore this curves gives rise to a (possibly non-trivial) additional element in the total holonomy group of
the (components of) $\mathcalp_1^{-1} (0)$ (resp. $\mathcalp_2^{-1} (0)$). In other words, unlike the preceding case, the total holonomy group of
the (components of) the singular fiber is no longer generated exclusively by local holonomy maps concentrated about singular points
of the mentioned fiber. This element is what was called an ``elliptic holonomy invariant''. In the present case, we have to keep track
of one or two of these invariants (according to the number of loops contained in the Dynkin diagram in question). We assume therefore that
these invariants coincide for the singular fibers $\mathcalp_1^{-1} (0), \, \mathcalp_2^{-1} (0)$. With this assumption in hand, the same argument employed
in Proposition~\ref{quasithm} applies to construct the desired (local) conjugacy between $\mathcalp_1, \, \mathcalp_2$.

It remains to consider the case in which $\mathcalp_1^{-1} (0), \, \mathcalp_2^{-1} (0)$ contain elliptic components. We denote by $E_1^1, E_1^2$ corresponding elliptic
components of $\mathcalp_1^{-1} (0), \, \mathcalp_2^{-1} (0)$. We note that these singular fibers may still contain a second elliptic component (which are going to be
denoted respectively $E_2^1, E_2^2$) or a loop of rational curve as above. More precisely, Dynkin diagrams containing
elliptic components are those labeled ``Type 1'', ``Type 12'', ``Type 13'' and ``Type 14'' in Ogg's classification \cite{ogg}. These elliptic components
may or or may not be pinched (representing a cusp). In any event an elliptic curve possesses two curves that are homotopically non-trivial so that it
gives rise to two ``elliptic holonomy invariants'' (clearly if the curve is pinched one of these two invariants is automatically trivial). So these invariants
also contribute to the total holonomy group of the singular fibers in question. However if they match for the corresponding singular fibers
$\mathcalp_1^{-1} (0), \, \mathcalp_2^{-1} (0)$ then the extension of the previous conjugacy is again well-defined. This completes the proof of the theorem.\qed

\vspace{0.1cm}

Let us close this section with the proof of the complement to Theorem~A. To do this, let us first show how the elliptic holonomy invariants can be
computed. Recall that as a divisor, the singular fiber is linearly
equivalent to a regular one. Set then $\mathcalp^{-1} (0) =\sum_{i=1}^k n_i C_i$ where $n_i \in \Z$. The coefficient $n_i$ is the {\it multiplicity}\, of
the corresponding component $C_i$. The computation of these coefficients offer no new difficult: they are already present in Ogg's
table \cite{ogg} and, in fact, can be determined by the same method exposed in Section~2.1. Therefore there is not loss of generality in resorting
to this information to compute the elliptic holonomy invariants.

To do this, note first that a genus~$2$ fibration has no multiple fiber, i.e. the greatest common divisor among all the $n_i$'s, $i=1 ,\ldots ,k$ is
always equal to~$1$. From this it follows that the multiplicity of a component $C_{i_0}$ of $\mathcalp^{-1} (0)$ is nothing but
the order of the ``total holonomy group of $C_{i_0}$
w.r.t. $\mathcalp$''. In fact, these numbers are nothing but the number of intersections of a nearby fiber with a
local transverse section $\Sigma$ through a point of $C_{i_0}$. We then conclude that the order of the holonomy maps arisen from ``elliptic holonomy invariants''
must divide the multiplicity of the component in question.

\vspace{0.1cm}

\noindent {\it Proof of the Complement to Theorem~A}. First we need to identify loops appearing in the Dynkin diagram associated to the singular fiber
in question. As mentioned it is the existence of these loops they may give rise to non-trivial elliptic holonomy invariants. Also these loops may be represented
by an (irreducible) elliptic curve appearing as a component of the singular fiber or by an actual loop of rational curves contained in the singular fiber.
A direct inspection on Ogg's list \cite{ogg}  shows that the multiplicity of an elliptic curve contained in a singular fiber can take only on the values~$1$ and~$2$.
So that the claim follows in this case.

A similar analysis applies to the cases where we have a loop of rational curves. In these case however a simple additional remark is needed.
First note that all singularities of the fibration appearing at singularities lying in a loop of rational curves are linear i.e. they have two eigenvalues
different from zero. In fact, it is enough to check Ogg's list again and recall that singularities are linear if and only if they are nodal singularities
for the (individual) analytic curve defined by the singular fiber itself. Now the possible elliptic holonomy invariant is necessarily associated
to the core path going through all these rational curves (the thread of the collar). In particular the elliptic holonomy invariant associated to this path
must divide the multiplicity of every rational curve constituting the collar in question. It is then immediate to check that all these elliptic holonomy invariants
are trivial ie. they equal~$1$. The statement is proved.\qed

As an example of application, consider the models of \cite{Ueno1} labeled as $[2 \! - \! 2I_{0-m}]$ (page 159) and $[3\!-\!II_{n-0} \ast]$ (page 172).
They exhibit isomorphic singular fibers containing an elliptic component of multiplicity~$2$. Thus the corresponding elliptic holonomy invariant
may be of order~$1$ (ie. the identity) or of order~$2$ (conjugate to a rotation of angle~$\pi$). These two possibilities account for the difference
between the corresponding fibrations containing the fibers in question. A completely analogous phenomenon occurs with the models
$[4\!-\!2I_n\!-\!m]$ and $[5\!-\!II_{n-p} \! \ast]$ appearing respectively in pages 181 and 184 of \cite{Ueno1}.

\section{Proofs for Theorems~B and~C}

In this last section $\mathcalp : M \rightarrow S$ (resp. $\mathcalp : M \rightarrow D \subset \C$) will stand for an isotrivial fibration of genus~$2$.
In the case we are dealing with a fibration over the disk $\mathcalp : M \rightarrow D \subset \C$, it is also supposed that $\mathcalp^{-1} (0)$
is the unique singular fiber of $\mathcalp$. The context of fibrations over the disk $D \subset \C$ is going to be referred to as a {\it semi-global
context (or fibration)}\, as opposed to a {\it global}\, fibration $\mathcalp : M \rightarrow S$ (or global context) where $S$ is a compact Riemann surface.

Let us first give the proof of Theorem~B. Fix isotrivial fibrations $\mathcalp_1 : M_1 \rightarrow D$, $\mathcalp_2 : M_2 \rightarrow D$ as in the
statement. We suppose that the typical fiber of $\mathcalp_1$ is isomorphic as Riemann surface to the typical fiber of
$\mathcalp_2$. Denote by $\sigma$ a holomorphic diffeomorphism between two such fibers. We are interested in the isotopy class of $\sigma$.
In this direction we have:

\begin{lema}
\label{FinishingLemma1.1}
The $\C^{\infty}$-diffeomorphism $H$ constructed in Theorem~A can be chosen so that its restriction from a regular fiber
of $\mathcalp_1$ to a regular fiber of $\mathcalp_2$ is isotopic to $\sigma$.
\end{lema}

\noindent {\it Proof}. Consider $\sigma$ as a diffeomorphism from a fiber $\mathcalp_1^{-1} (x_1)$ to a fiber $\mathcalp_2^{-1} (x_2)$.
Consider also the line $l_1$ joining $x_1$ to $0 \in D\subset \C$ (resp. $l_2$ joining $x_2$ to $0 \in D\subset \C$). By using
the fact that $\mathcalp_1, \mathcalp_2$ are $C^{\infty}$-fibrations away from their respective singular fibers, the isotopy class of $\sigma$
can be ``transported over $l_1, l_2$'' to induce a diffeomorphism between the singular fibers $\mathcalp_1^{-1} (0)$ and
$\mathcalp_2^{-1} (0)$. Indeed regular fibers of $\mathcalp_1$ (resp. $\mathcalp_2$) are ramified coverings of $\mathcalp_1^{-1} (0)$
(resp. $\mathcalp_2^{-1} (0)$) with fixed degree and ``bounded ramification''.
Though this diffeomorphism is not canonically defined, its isotopy class is so. We consider then a diffeomorphism
$\textsf{h}$ between $\mathcalp_1^{-1} (0), \, \mathcalp_2^{-1} (0)$ whose isotopy class agrees, in the above sense, with the isotopy class
of $\sigma$.

To conclude the proof of the lemma, it suffices to construct the diffeomorphism $H$ of Theorem~A by considering the identification
of the singular fibers $\mathcalp_1^{-1} (0), \, \mathcalp_2^{-1} (0)$ given by $\textsf{h}$. In fact, the method used in Theorem~A depends on
the lifts of paths in these singular fibers and the correspondence between these paths is settled by $\textsf{h}$. It then becomes
clear that $H$ satisfies the required condition.\qed

We also have the following well-known lemma:

\begin{lema}
\label{FinishingLemma2.2}
Let $S$ denote a compact hyperbolic Riemann surface. Suppose that $\sigma_1, \sigma_2$ are two automorphisms of $S$
whose induced actions on $H^1 (S, \C)$ turns out to coincide. Then $\sigma_1 =\sigma_2$. In particular two isotopic automorphisms of
$S$ must coincide.
\end{lema}

\noindent {\it Proof}. It suffices to prove that an automorphism $\sigma$ of $S$ acting trivially on the cohomology ring of $S$
must coincide with the identity.
By using the hyperbolic structure of $S$ we can define a Jacobian map ${\rm Jac}$ from $S$ to the
Albanese torus $\mathbb{T}$ of $S$ which is equivariant with respect to the group $G$ of automorphisms of $S$. 
More precisely $G$ acts on $S$ in the obvious way and on $\mathbb{T}$ by affine automorphisms. In particular,
since $\sigma$ acts trivially on the cohomology of $S$, its action on $\mathbb{T}$ coincides with a translation. Therefore
either this translation is trivial or it has no fixed point. In the latter case the equivariance of ${\rm Jac}$ implies that
$\sigma$ has no fixed points in $S$ as well. This is however impossible since the Lefschetz number of $\sigma$ is strictly
positive.\qed

\vspace{0.1cm}

\noindent {\it Proof of Theorem~B}. Let $\mathcalp_1 : M_1 \rightarrow D$ and $\mathcalp_2 : M_2 \rightarrow D$ be as in the statement. According to
Theorem~A and Lemma~\ref{FinishingLemma1.1} there exists a $\C^{\infty}$-diffeomorphism $H$ conjugating the fibrations
$\mathcalp_1, \, \mathcalp_2$. Furthermore $H$ is transversely holomorphic and the $C^{\infty}$ diffeomorphism induced by
$H$ between corresponding regular fibers $\mathcalp_1^{-1} (x_1), \, \mathcalp_2^{-1} (x_2)$ of $\mathcalp_1, \, \mathcalp_2$ is isotopic to a holomorphic
diffeomorphism from $\mathcalp_1^{-1} (x_1)$ to $\mathcalp_2^{-1} (x_2)$.

To construct a holomorphic conjugacy between $\mathcalp_1, \, \mathcalp_2$ we proceed as follows. Consider the circle parametrized by
$x_1 e^{2\pi it}$, $t \in [0,1]$. We shall produce a fibered deformation of $H$ over this circle so that the resulting diffeomorphism
will induce holomorphic diffeomorphisms between the corresponding fibers. To do this, we begin by deforming $H$ through a $C^{\infty}$-isotopy
on the fiber $\mathcalp_1^{-1} (x_1)$ so that it becomes a holomorphic diffeomorphism from $\mathcalp_1^{-1} (x_1)$ to $\mathcalp_2^{-1} (x_2)$. Let
$\sigma$ denote this holomorphic diffeomorphism. We then continue deforming $H$ on the fibers of $\mathcalp_1$ lying over the circle
$x_1 e^{2\pi it}$, $t \in [0,1]$ so that it induces holomorphic diffeomorphisms between corresponding fibers. We need to show that this
deformation is well-defined after winding around the origin. For this let $\tilde{\sigma}$ denote the holomorphic diffeomorphism from 
$\mathcalp_1^{-1} (x_1)$ to $\mathcalp_2^{-1} (x_2)$ induced after one turn over the mentioned circle. We need to show that $\tilde{\sigma}= \sigma$. This however
follows from Lemma~\ref{FinishingLemma2.2}: since the initial $C^{\infty}$-diffeomorphism $H$ is globally defined, it is clear that both
$\tilde{\sigma}, \, \sigma$ are isotopic to the restriction of $H$ to $\mathcalp_1^{-1} (x_1)$. Thus $\tilde{\sigma}$ is itself isotopic to
$\sigma$ and hence $\tilde{\sigma}= \sigma$ as desired.

With the above procedure, we have construct a $C^{\infty}$-diffeomorphism $\tilde{H}$ conjugating $\mathcalp_1, \, \mathcalp_2$
satisfying the following conditions:
\begin{enumerate}
\item $\tilde{H}$ is transversely holomorphic and induces holomorphic diffeomorphisms between the regular fibers of
$\mathcalp_1, \, \mathcalp_2$.

\item $\tilde{H}$ admits a continuous extension to the singular fiber $\mathcalp_1^{-1} (0)$.
\end{enumerate}
The assertion that $\tilde{H}$ admits a continuous extension to $\mathcalp_1^{-1} (0)$ is itself a consequence of the proof of
Theorem~A. In fact $H$ was constructed by means of a fibration ``transverse'' to the singular fiber which can be used to
compare the fiberwise isotopies for fibers near the singular one. Now the classical Hartogs theorem ensures that $\tilde{H}$
is holomorphic on $M_1 \setminus \mathcalp_1^{-1} (0)$. In turn Riemann extension allows us to conclude that $\tilde{H}$
is holomorphic on $\mathcalp_1^{-1} (0)$ as well. The theorem is proved.\qed

In the remainder of this section we are going to discuss the structure of (global) isotrivial fibrations. Our discussion relies
heavily on the following observation applying to both semi-global and global contexts.

\begin{prop}
\label{Lastsection1.1}
Let $\mathcalp : M \rightarrow S$ (resp. $\mathcalp : M \rightarrow D \subset \C$) be an isotrivial fibration of genus~\textsf{g} and denote by
$\mathcalp^{-1} (p_i)$ its singular fibers, $i=1, \ldots , k$. Then there exists a
holomorphic foliation $\mcd$ defined on $M$ and satisfying the following conditions:
\begin{enumerate}
\item $\mcd$ leaves invariant the singular fibers $\mathcalp^{-1} (p_i)$, $i=1, \ldots , k$, of $\mathcalp$.

\item Away from $\bigcup_{i=1}^k \mathcalp^{-1} (p_i)$, the foliation $\mcd$ is transverse to $\mathcalp$.

\end{enumerate}
\end{prop}

\noindent {\it Proof}. We shall first construct $\mcd$ on $M \setminus \bigcup_{i=1}^k \mathcalp^{-1} (p_i)$. For this recall that every
point in $S \setminus \{p_1, \ldots ,p_k\}$ is contained in a neighborhood $W \subset S$ for which there exists a local trivialization
$$
\mathcalp^{-1} (W) \simeq W \times S_{\textsf{g}}
$$
where $S_{\textsf{g}}$ stands for the Riemann surface determined by the regular fibers of $\mathcalp$. The existence of this trivialization
is a consequence of a well-known theorem due to Fischer and Grauert (see \cite{bpv}). In the trivializing coordinates, we consider the
foliation defined by horizontal lines. We need to check that these local foliations patch together in a
foliation defined on $M \setminus \bigcup_{i=1}^k \mathcalp^{-1} (p_i)$.
For this let us consider two trivializing coordinates as above over open sets $W_1, W_2$ contained in $S$. Let $z$ be a local
coordinate defined on the intersection $W_1 \cap W_2$ the change of trivializations gives rise to a diffemorphism $\varphi :
W_1 \cap W_2 \times S_{\textsf{g}} \rightarrow W_1 \cap W_2 \times S_{\textsf{g}}$ of the form
$$
\varphi (z,p) = (\varphi_1 (z) , \xi (z,p)) \, .
$$
In particular for $z$ fixed the induced map $p \mapsto  \xi (z,p))$ is an automorphism of the Riemann surface $S_{\textsf{g}}$. The horizontal
vector field $\partial /\partial z$ is pulled back to the vector field
$$
\frac{1}{\varphi_1' (z)} \frac{\partial}{\partial z} - \frac{\partial \xi /\partial z}{\varphi_1'(z) \partial \xi /\partial y} \frac{\partial}{\partial y}
$$
where $y$ is a local coordinate about $p \in S_{\textsf{g}}$. However the group of automorphism of $S_{\textsf{g}}$ is finite, and therefore
discrete, since $\textsf{g} \geq 2$. Because $\xi (z,p)$ can be regarded as a family of automorphisms of $S_{\textsf{g}}$ varying continuously with
$z$, it follows that this family must be constant. Hence the partial derivative $\partial \xi /\partial z$ vanishes identically and we conclude that the
above vector field is still horizontal. The horizontal direction is therefore well-defined and endows a holomorphic foliation defined on
 $M \setminus \bigcup_{i=1}^k \mathcalp^{-1} (p_i)$.
 
To finish the proof of the proposition we only need to prove that the above constructed foliation can continuously be extended to the regular
part of the singular fibers of $\mathcalp$ and that this extension is tangent to the mentioned singular fibers. Indeed in this case Riemann extension 
applies to ensure the foliation is actually holomorphic away from the singular points of the singular fibers of $\mathcalp$. The latter set being of
codimension~$2$, it follows that the foliation is holomorphic on all of $M$.

To check the claim, consider a point $P \in \mathcalp^{-1} (p_i) \subset M$ that is regular for the singular fiber in question. Since $\mathcalp$ is not
a submersion at $P$, it promptly follows that the tangent spaces to the leaves of the above constructed foliation converge (in direction)
towards the tangent space of $\mathcalp^{-1} (p_i)$ at $P$. Thus the foliation can continuously be extended to a neighborhood of $P$ by saying
that its local leaf through $P$ coincides with a neighborhood of $P$ in $\mathcalp^{-1} (p_i)$. The proposition is proved.\qed

By applying Proposition~\ref{Lastsection1.1} to a semi-global fibration $\mathcalp : M \rightarrow D \subset \C$, we obtain the
following:

\begin{coro}
\label{Lastsection2.2}
Suppose that $\mathcalp : M \rightarrow D \subset \C$ is a semi-global isotrivial fibration of genus $\textsf{g} \geq 2$. Then the
monodromy of $\mathcalp$ about its singular fiber can be represented by a holomorphic diffeomorphism of the typical fiber
$S_{\textsf{g}}$.
\end{coro}

\noindent {\it Proof}. Consider the foliation $\mcd$ constructed above and fix $z_0 \in D$, $z_0 \neq 0$. The circle $t \mapsto e^{it} z_0$,
$t \in [0, 2\pi]$ can be lifted in the leaves of $\mcd$ with respect to the fibration $\mathcalp$. Because $\mcd$ is a holomorphic
foliation, this lifting gives rise to a holomorphic diffeomorphism of $\mathcalp^{-1} (z_0)$ (the holonomy induced by the fundamental group of
$D^{\ast}$). Clearly this diffeomorphism is a representative of the monodromy of $\mathcalp$.\qed

More generally let $\mathcalp : M \rightarrow S$ denote an isotrivial fibration of genus~$2$ whose critical values are $\{p_1, \ldots ,p_s \} \subset S$.
The fundamental group of the punctured surface $S \setminus \{p_1, \ldots ,p_s \}$ is generated by a generating set for the fundamental group
of $S$ plus small loops around each of the critical values $p_i$, $i=1, \ldots ,s$. This group will be denoted by $\pi_1 (S \setminus \{p_1, \ldots ,p_s \})$.
The existence of the foliation $\mcd$ gives rise to a representation from $\pi_1 (S \setminus \{p_1, \ldots ,p_s \})$ in the group ${\rm Aut}\, (S_{\textsf{g}})$
of automorphism of the typical fiber of $\mathcalp$. The element of ${\rm Aut}\, (S_{\textsf{g}})$ produced by the holonomy corresponding to a small
circle around a critical value is nothing but the holomorphic representative of the monodromy of the singular fiber in question (cf. Corollary~\ref{Lastsection2.2}).
To compute the (necessarily finite) order of this element is easy as shown by the next lemma.

\begin{lema}
\label{Lastsection3.3}
Let $\mathcalp : M \rightarrow D \subset \C$ be as in Corollary~\ref{Lastsection2.2} and denote by $\sigma$ the element of ${\rm Aut}\, (S_2)$
generated by the monodromy of $\mathcalp^{-1} (0)$. If $\mathcalp^{-1} (0)$ has only normal crossings, then the order of $\sigma$ is the least common
multiple ${\rm l.c.m.}$ of the multiplicities $n_i$ of the irreducible components of $\mathcalp^{-1} (0)$. In general, modulo performing finitely many
blow-ups (cf. Section~2.2 and Proposition~\ref{localmodels}) $\mathcalp^{-1} (0)$ can be turned into a (non-minimal) fiber possessing only
normal crossings.
\end{lema}

\noindent {\it Proof}. A consequence of Lemma~\ref{FinishingLemma2.2} together with Corollary~\ref{Lastsection2.2} is that, as far as isotrivial
fibrations are concerned, there is no distinction between standard monodromy and ``non-commutative monodromy'' in the sense of
\cite{MMontesinos}. In particular if the standard monodromy is trivial so is the ``non-commutative monodromy'' and, in this case,
$\mathcalp^{-1} (0)$ is a regular fiber.

Now suppose that $\mathcalp^{-1} (0)$ has only normal crossing singular points.
Let $n$ be the least common multiple of the elements in the set formed by the multiplicities $n_i$ of the irreducible components of $\mathcalp^{-1} (0)$.
Next let $\mathcalp^{(n)} : M^{(n)} \rightarrow D$ be the $n$-root fibration of $\mathcalp$, cf. \cite{bpv}. Then $\mathcalp^{(n)}$ has a stable fiber sitting over
$0 \in \C$. Again the monodromy associated to this stable fiber must be finite so that $\mathcalp^{(n)}$ is a root of a stable fibration whose monodromy is trivial.
It then follows that the monodromy of $\mathcalp^{(n)}$ itself must be trivial (in particular $\mathcalp^{(n)}$ is not singular). The rest of the argument is obvious.\qed

The standard stable reduction theorem combined to what precedes yields part of the statement of Theorem~A. As already mentioned modulo performing
finitely many blow-ups, we can always suppose that all singular fibers have at worst normal crossings.

\begin{coro}
\label{Lastsection4.4}
Let $\mathcalp : M \rightarrow S$ be an isotrivial fibration of genus~$2$ all of whose singular fibers have only normal crossing singularities. Then there exists
a cyclic covering $\widetilde{S} \rightarrow S$, ramified only over the critical values of $\mathcalp$ in $S$ and over a single additional point in $S$ such that
the fiber product $\widetilde{S} \times_S M$ is bimeromorphically equivalent to a regular fibration of genus~$2$ (i.e. a Kodaira fibration, necessarily isotrivial).
\end{coro}

It is now easy to rule out some models of singular fibers that cannot occur in an isotrivial genus~$2$ fibration. In fact, further information can be obtained
if the group of autormorphism ${\rm Aut}\, (S_2)$ is brought to bear. Some information on this group is collected below, the reader may check
\cite{accola} for further details and proofs. First it is well-known that the hyperelliptic involution $\jmath$ is {\it unique}\, and central in ${\rm Aut}\, (S_2)$.
The quotient group ${\rm Aut}\, (S_2) / \langle \jmath \rangle$ is a finite group of automorphism of the Riemann sphere $\overline{\C}$. The latter groups
are classified and the only possibilities are the cyclic groups, the dihedral group $D_n$, the alternating groups $\mathcal{A}_4$ and $\mathcal{A}_5$
and the symmetric group $\mathcal{S}_4$.

\vspace{0.2cm}

\noindent {\it Proof of Theorem~C}. First recall that the order of ${\rm Aut}\, (S_2)$ belongs to the interval $[24, 84]$ (the upper bound being the usual
Hurwitz bound).  This immediately allows us to restricted the above mentioned list of groups to cyclic groups, $D_n$ for $n=6, \ldots ,21$, $\mathcal{A}_4$
and $\mathcal{S}_4$. It also well-known that the order of an automorphism of $S_2$ cannot exceed~$5$ provided that this order is prime. The classical
Sylow theorem then implies that the order of ${\rm Aut}\, (S_2)$ can have only~$2, 3$ and~$5$ as prime factors. Thus we conclude that
${\rm Aut}\, (S_2)$ is as indicated in the statement.

Finally recall that we have a representation from $\pi_1 (S \setminus \{p_1, \ldots ,p_s \})$ in the group ${\rm Aut}\, (S_{\textsf{g}})$, where
$\{p_1, \ldots ,p_s \} \subset S$ stands for the critical values of $\mathcalp$. In particular the monodromy generated by a loop around $p_j$
belongs to the image of this representation and thus is an element of ${\rm Aut}\, (S_{\textsf{g}})$. Thanks to Lemma~\ref{Lastsection3.3},
we know that the order of this element is nothing but the least common multiple of the multiplicities $n_j^i$ of the irreducible components
of $\mathcalp^{-1} (p_j)$. Hence all these multiplicities  cannot have prime factors different from~$2, \, 3$ or~$5$. Besides if there is a component
having~$5$ as prime factor, the order of ${\rm Aut}\, (S_{\textsf{g}})$ must be divisible by~$5$ as well. The theorem is proved.\qed

\vspace{0.1cm}

\noindent {\bf A comment on the notion of dual fibers}: We may define two singular fibers to be {\it dual}\, if there exists a genus~$2$
fibration $\mathcalp : M \rightarrow \overline{\C}$ possessing exactly two singular fibers that are respectively isomorphic to the initial ones.
Here $\overline{\C}$ stands for the Riemann sphere. Unlike the case of {\it elliptic fibrations}\, all genus~$2$ fibrations over
$\overline{\C}$ possessing only two singular fibers are automatically isotrivial. The proof is well-known: by choosing a suitable
symplectic basis, we can construct a morphism from $\C$ (viewed as the universal covering of $\overline{\C}$ minus two points)
in the Siegel half-space $\mathcal{H}_2$. If the fibration were not isotrivial, the classical Torelli theorem would imply that this
morphism is not constant. This however contradicts Liouville theorem since $\mathcal{H}_2$ is isomorphic to a bounded domain.

In particular the arguments presented above can be applied to define a match between singular fibers admitting duals. We also note
that many examples of these situations can directly be obtained from the ``local models'' given in Lemma~\ref{Section2.11} and
in Proposition~\ref{Section2.22} along the same lines of the example constructed at the end of Section~2.2 with the local
model $x^5y^4(x-y)$.


\subsection*{Acknowledgment}
We are very much indebted to the anonymous referee of this article who provided us with a number of references
and valuable information on fibrations of low genus, including open questions and general strategies. In particular this material allowed us to
improve on a preliminary version and to give a more reasonable perspective of the topic concerning fibrations of low genus.

We are also grateful to J.-F. Mattei for discussions concerning the variation of the
complex structure of fibers and to G. Tian for the interest showed in this paper.

\end{document}